\input amstex
\magnification=\magstep1
\documentstyle{amsppt}
\advance\hoffset by 0.1 truein
\TagsOnRight
\hoffset=0.29in
\hsize=4.75in
\voffset=0in
\vsize=7.5in
\parindent=10pt

\font\bigrm=cmr12
\font\smallrm=cmr8
\document

\  

\vskip0.8cm

\centerline{\bigrm SEPARATION PROPERTIES OF THETA FUNCTIONS}

\vskip0.5cm

\centerline{\rm EDUARDO ESTEVES\footnote 
{Research supported by an MIT Japan Program 
Starr fellowship, and by CNPq, Proc. 300004/95-8}}

\vskip0.5cm

\author 
\smallrm EDUARDO ESTEVES
\endauthor

\title 
\smallrm SEPARATION PROPERTIES OF THETA FUNCTIONS
\endtitle

\vskip0.5cm

\bf 0.\hskip0.2cm Introduction.\hskip0.4cm \rm Let $X$ be a 
non-singular, connected, projective curve defined over an 
algebraically closed field $k$. Let $U^s$ denote the set of isomorphism 
classes of stable vector bundles on $X$ with given degree $d$ and 
rank $r$. In the sixties, C.S. Seshadri and D. Mumford (\cite{\bf 14}, 
\cite{\bf 17} and \cite{\bf 18}) supplied $U^s$ with a natural 
structure of quasi-projective variety, together with a natural 
compactification $U$, by adding semistable vector bundles at the boundary. 
The method used in the construction of such a structure was Mumford's then 
recently developed Geometric Invariant Theory \cite{\bf 8}. 
Roughly, the method 
consists of producing a variety $R$, and an action of a reductive group 
$G$ on $X$, linearized at some ample invertible sheaf $L$ on $R$, such that 
$U=R/G$ set-theoretically. Then, Geometric Invariant Theory (G.I.T.) 
tells us how to supply $R/G$ with 
a natural scheme structure, obtained from the $G$-invariant sections of 
tensor powers of $L$.

Up until recently, Seshadri's and Mumford's construction was 
the only purely algebraic construction available. In 1993, Faltings 
\cite{\bf 7} showed how to construct $U^s$, and 
its compactification $U$, avoiding G.I.T.. His method, described also in 
\cite{\bf 20}, consisted in 
considering the so-called theta functions on $R$, 
naturally defined provided $R$ admits a family with 
the so-called local universal property. (We observe that the theta functions 
considered in this article are just those associated with vector 
bundles on $X$, as it will be clear from our definition in Sect. 2. 
Beauville \cite{\bf 3, \rm Sect. 2} has a more encompassing definition of 
theta functions than ours.) The theta functions are 
in fact $G$-invariant sections of tensor powers of a certain $G$-linear 
invertible sheaf $L'_{\theta}$ on $R$. Roughly speaking, using 
his first main lemma \cite{\bf 20\rm , Lemma 3.1, p. 166}, 
Faltings showed that there are enough theta functions 
to produce a $G$-invariant morphism, $\theta\: R @>>> \text{\bf P}^N$. 
By semistable reduction, the image,
$$
U_{\theta}:=\theta(R)\subseteq\text{\bf P}^N,
$$
is a closed 
subvariety. Since $\theta$ is $G$-invariant, then $\theta$ factors 
(set-theoretically) through a map, $\pi\: U=R/G @>>> U_{\theta}$. Then, 
using his 
second main lemma \cite{\bf 20\rm , Lemma 4.2, p. 174}, 
Faltings showed that there is 
a bijection between the normalization of $U_{\theta}$ and $U$, through which 
$U$ acquires a natural structure of projective variety.

There are interesting consequences of Faltings' work; for instance, we can 
get from his first main lemma a cohomological characterization of 
semistable bundles (cf. Thm. 2 or \cite{\bf 20\rm , Thm. 6.2, p. 187}).

In his description of Faltings' construction, Seshadri posed the 
following question (\cite{\bf 20\rm , Rmk. 6.1, p. 187}): 
Let $L_{\theta}$ denote the ample sheaf on $U$ 
lying below $L'_{\theta}$. The theta functions of powers of $L'_{\theta}$ 
descend to sections of powers of $L_{\theta}$, and span a graded subalgebra 
$A_{\theta}$ of
$$
B_{\theta}:=\bigoplus_{m\geq 0} \Gamma(U,L_{\theta}^{\otimes m}).
$$
How close is $A_{\theta}$ to $B_{\theta}$? Lest the latter question is too 
difficult, we can think more 
geometrically, and rephrase it: 
how close is the normalization map,
$$
\pi\: U @>>> U_{\theta},
$$
to being an isomorphism, if we consider all 
theta functions? In other words, how much of the moduli space $U$ can theta 
functions describe? This question arises naturally from Faltings' 
construction, and is relevant: if $U=U_{\theta}$, then we would have a 
better handle on a canonical projective embedding of $U$. 

The goal of the present article is to provide a partial answer to 
Seshadri's question. 
Our main technical lemma (Lemma 4) is a generalization of Faltings' 
first main lemma. 
From Lemma 4, we obtain a cohomological 
characterization of stability (Thm. 6), allowing us to obtain a 
quick proof of the fact that $U$ is a fine moduli space 
if the degree $d$ and rank $r$ are coprime (Cor. 8). Our proof of this fact 
avoids completely G.I.T.. 
Using Lemma 4, we obtain separation lemmas (Lemma 9 and Lemma 11), 
allowing us to show that $\pi$ is bijective (Thm. 15) and 
(if the characteristic of the ground field $k$ is 0) an isomorphism over 
$U^s$ (Thm. 18).

Actually, motivated by the recent interest in the 
compactification of relative Jacobians over families of singular curves 
(\cite{\bf 4}, \cite{\bf 16} and the more general \cite{\bf 21}), 
and to provide support to 
\cite{\bf 6}, 
I have tried to be as general as possible. So, all the results of the 
present paper (with the exception of Thm. 18) apply to 
a general projective, connected, reduced curve $X$, with $n$ irreducible 
components, defined over an algebraically closed field $k$. 
In \cite{\bf 19\rm , Part 7}, 
Seshadri constructed moduli 
spaces $U(\underline a,\underline m,\chi)$ of 
torsion-free sheaves with multirank $\underline m=(m_1,\dots,m_n)$ and 
Euler characteristic $\chi$ that are semistable with respect to a 
polarization $\underline a=(a_1,\dots,a_n)$. 
In this article, we consider theta functions on 
$U(\underline a,\underline m,\chi)$ and ask ``Seshadri's question'': 
how close is the rational map,
$$
\pi\: U(\underline a,\underline m,\chi) @>>> 
U_{\theta}(\underline a,\underline m,\chi),
$$
obtained from all theta functions, to being an isomorphism? If 
$n=1$ and $m=1$, then $\pi$ has already been shown to be an 
isomorphism in \cite{\bf 5}. In general, our best results to this 
extent are Thm. 15 and Thm. 16. In particular, if 
$\underline m=(1,\dots,1)$ (the compactified Jacobian case), 
then we obtain that $\pi$ is bijective and an isomorphism over 
the stable locus, $U^s(\underline a,\underline m,\chi)\subseteq 
U(\underline a,\underline m,\chi)$.

\vskip0.8cm

\bf 1.\hskip0.2cm Preliminaries.\hskip0.4cm \rm Let $X$ be a curve, that is, 
a projective, connected, reduced scheme of pure dimension 1 over 
an algebraically closed field $k$. Let $X_1,\dots,X_n$ denote the 
irreducible components of $X$. All schemes are assumed to be 
locally of finite type over $k$. 
By a point we mean a closed point. By a vector bundle we 
mean a locally free sheaf of constant rank.

If $H$ is a coherent sheaf on $X$, we let $\underline r_H$ denote its 
\sl multirank\rm ; more precisely, 
$\underline r_H:=(r_1(H),\dots,r_n(H))$, where 
$r_i(H)$ is the generic rank of $H$ on $X_i$ for every $i$. We let 
$r_H$ denote the maximum generic rank of $H$. 
If $E$ is a vector bundle on $X$, we let $\underline d_E$ denote its \sl 
multidegree\rm , that is, $\underline d_E:=(\deg_{X_1} E,\dots,\deg_{X_n} E)$. 
If $H$ is a coherent sheaf, and $E$ is a vector bundle on $X$, then
$$
\chi(E\otimes H)=r_E\chi(H)+\underline d_E\cdot\underline r_H.
$$
In particular, note that $\chi(E\otimes H)$ depends only on the rank and 
multidegree of $E$, rather than on $E$ itself.

A \sl torsion-free \rm sheaf on $X$ is a coherent sheaf with 
no embedded components. The dualizing sheaf of $X$, denoted $\omega$, is 
torsion-free \cite{\bf 1\rm , (6.5), p. 95}. Given any 
coherent sheaf $H$, we let $H^{\omega}:=
\underline{\text{Hom}}_X(H,\omega)$. 
Since $\omega$ is torsion-free, then so is $H^{\omega}$. There is a 
natural surjective homomorphism $H \twoheadrightarrow H^{\omega\omega}$, whose 
kernel is the torsion subsheaf of $H$. 

Let $S$ be a scheme. A coherent $S$-flat sheaf $\Cal I$ on $X\times S$ 
is called 
\sl relatively torsion-free over $S$ \rm if the fibre $\Cal I(s)$ is 
torsion-free for every $s\in S$. Given any coherent 
sheaf $\Cal H$ on $X\times S$, we let
$$
\Cal H^{\omega}:=\underline{\text{Hom}}_{X\times S}(\Cal H,\omega\otimes
\Cal O_S).
$$
If $\Cal I$ is relatively torsion-free on $X\times S$ over $S$, then 
it follows from \cite{\bf 1\rm , (1.9), p. 59} and 
\cite{\bf 1\rm , (6.5), p. 95} that $\Cal I^{\omega}$ is also 
relatively torsion-free, with $\Cal I^{\omega}(s)=\Cal I(s)^{\omega}$ for 
every $s\in S$. Anyhow, if $\Cal H$ is a coherent sheaf on 
$X\times S$, it follows from loc. cit. 
that there is an open dense subscheme $S'\subseteq S$ such that 
$\left.\Cal H^{\omega}\right|_{X\times S'}$ is relatively torsion-free with 
$\Cal H^{\omega}(s)=\Cal H(s)^{\omega}$ for every $s\in S'$.

If $L$ is an ample invertible sheaf on $X$, we let 
$P_H:=P_H(T)$ denote the Hilbert polynomial 
of a coherent sheaf $H$ on $X$ with respect to $L$. 
Using the dualizing properties of $\omega$, 
we get that $P_H(T)=-P_{H^{\omega}}(-T)$ for 
every coherent sheaf $H$ on $X$.

\vskip0.8cm

\bf 2.\hskip0.2cm Theta functions.\hskip0.4cm \rm Let $S$ be a scheme. 
Let $\Cal F$ be a coherent sheaf on $X\times S$ that is flat over $S$. 
The \sl determinant of cohomology of $\Cal F$ over $S$ \rm 
is the invertible sheaf $\Cal D(\Cal F)$ on
$S$ constructed as follows: locally on $S$ there is a
complex,
$$
0 @>>> G^0 @>\lambda>> G^1 @>>> 0,
$$
of free sheaves of finite rank such that, for every coherent sheaf $M$
on $S$, the cohomology groups of $G^{\bullet}\otimes M$
are equal to the higher direct images of $\Cal F\otimes M$
under the projection $p\: X\times S @>>> S$. 
The complex $G^{\bullet}$ is unique (up to
unique quasi-isomorphism). Hence, its determinant,
$$
\det G^{\bullet}:=(\bigwedge^{\text{rank } G^1}
G^1)\otimes(\bigwedge^{\text{rank } G^0} G^0)^{-1},
$$
is unique (up to canonical isomorphism). 
The uniqueness allows us to glue together the
local determinants to obtain the invertible sheaf $\Cal
D(\Cal F)$ on $S$. 

The most important properties of the determinant of cohomology are:

(a)\hskip0.2cm \sl Additive property\rm : If
$$
\alpha\: 0 @>>> \Cal F_1 @>>> \Cal F_2 @>>> \Cal F_3 @>>> 0
$$
is a short exact sequence of $S$-flat coherent sheaves on $X\times S$, 
then there is a naturally associated isomorphism:
$$
\Cal D_{\alpha}\: 
\Cal D(\Cal F_2) \cong \Cal D(\Cal F_1)\otimes\Cal D(\Cal F_3).
$$

(b)\hskip0.2cm \sl Projection property\rm : If 
$L$ is a line bundle on $S$, and 
$\chi(\Cal F(s))=d$ for every $s\in S$, then there is a naturally 
associated isomorphism:
$$
\Cal D_L\: \Cal D(\Cal F\otimes L)\cong\Cal D(\Cal F)\otimes L^{\otimes d}.
$$

(c)\hskip0.2cm \sl Base-change property\rm : If 
$\nu\: T @>>> S$ is any morphism, 
then there is a naturally associated base-change isomorphism:
$$
\Cal D_{\nu}\: \Cal D((\text{id}_X,\nu)^*\Cal F)\cong\nu^*\Cal D(\Cal F).
$$

For a more systematic development of the theory of
determinants, see \cite{\bf 11}. It is also possible to adopt
a more concrete approach to define $\Cal D(\Cal F)$, like the one 
used in \cite{\bf 2\rm , Ch. IV, \S 3}. 

If $\Cal F$ is an $S$-flat coherent sheaf on $X\times S$ with 
$\chi(\Cal F(s))=0$ for every $s\in S$, then 
there is a canonical global section $\sigma_{\Cal F}$ of $\Cal
D(\Cal F)$ that is constructed as follows: since 
$\chi(\Cal F(s))=0$ for every $s\in S$, then the ranks of
$G^0$ and $G^1$ in the local complex $G^{\bullet}$ are
equal. By taking the determinant of $\lambda$ we obtain a 
section of $\det G^{\bullet}$. Since the complex $G^{\bullet}$ is unique, 
such section is also unique, allowing us to glue the
local sections to obtain $\sigma_{\Cal F}$. 
The zero locus of $\sigma_{\Cal F}$ on $S$
parametrizes the points $s\in S$ such that
$$
h^0(X,\Cal F(s))=h^1(X,\Cal F(s))=0.
$$
(Another way of viewing $\sigma_{\Cal F}$ is as a generator of the 
0-th Fitting ideal of $R^1p_*\Cal F$.)

The global section $\sigma_{\Cal F}$ satisfies properties 
compatible with those of $\Cal D(\Cal F)$. For instance, 
we have the 
additive property: if 
$$
\alpha\: 0 @>>> \Cal F_1 @>>> \Cal F_2 @>>> \Cal F_3 @>>> 0
$$
is a short exact sequence of $S$-flat coherent sheaves on $X\times S$ of 
relative Euler characteristic $0$ over $S$, then
$$
\sigma_{\Cal F_2}=\sigma_{\Cal F_1}\otimes\sigma_{\Cal F_3}
$$
under the identification given by $\Cal D_{\alpha}$. 
We leave it to the reader to state 
the projection and base-change properties of the global sections 
$\sigma_{\Cal F}$.

Let $E$ be a vector bundle on $X$. Let 
$S$ be a scheme, and $\Cal H$ be an $S$-flat coherent sheaf on 
$X\times S$. Assume that $\chi(E\otimes\Cal H(s))=0$ for every $s\in S$. 
We define:
$$
\Cal L_E(\Cal H):=\Cal D(E\otimes\Cal H) \text{\  \  and \  \  } 
\theta_E(\Cal H):=\sigma_{E\otimes\Cal H}.
$$
The line bundle $\Cal L_E(\Cal H)$ is called a \sl theta line bundle\rm , 
and $\theta_E(\Cal H)$ is called its \sl theta function\rm . 
We let $\Theta_E(\Cal H)\subseteq S$ 
denote the zero-scheme of $\theta_E(\Cal H)$, and call it a 
\sl theta divisor\rm . 

\vskip0.4cm

\smc Lemma 1. \rm (Faltings)\hskip0.4cm \sl Let $S$ be a scheme. Let 
$\Cal I$ and $\Cal J$ be $S$-flat coherent sheaves on $X\times S$. 
Let $E$ and $F$ be vector bundles on $X$. Assume that:

{\rm (a)}\hskip0.2cm $\chi(\Cal I(s))=\chi(\Cal J(s))$ and 
$\underline r_{\Cal I(s)}=\underline r_{\Cal J(s)}$ for 
every $s\in S$;

{\rm (b)}\hskip0.2cm $r_E=r_F$ and $\det E\cong\det F$.

\parindent=0pt

Then, there is a canonical isomorphism,
$$
\Phi_{F,E}\: \Cal D(\Cal I\otimes F)\otimes\Cal D(\Cal J\otimes E)\cong
\Cal D(\Cal I\otimes E)\otimes\Cal D(\Cal J\otimes F),
$$
whose formation commutes naturally with base change. In 
addition, $\Phi_{F,E}$ is additive on $E$ and $F$, in the following sense: 
If
$$
\alpha\: 0 @>>> E_1 @>>> E_2 @>>> E_3 @>>> 0 \text{\  \  and \  \  }
\beta\: 0 @>>> F_1 @>>> F_2 @>>> F_3 @>>> 0
$$
are short exact sequences of vector bundles on $X$ such that 
$r_{E_i}=r_{F_i}$ and $\det E_i\cong\det F_i$ for $i=1,2,3$, then
$$
(\Phi_{F_1,E_1}\otimes\Phi_{F_3,E_3})\circ
(\Cal D_{\beta\otimes\Cal I}\otimes\Cal D_{\alpha\otimes\Cal J})=
(\Cal D_{\alpha\otimes\Cal I}\otimes\Cal D_{\beta\otimes\Cal J})\circ
\Phi_{F_2,E_2}.
$$

\vskip0.4cm

\parindent=10pt

\sl Proof.\hskip0.4cm \rm It follows from (a) that 
the sheaves $\Cal I$ and $\Cal J$ are 
generically isomorphic, hence coincide generically in $K$-theory, in the 
sense defined by Faltings in \cite{\bf 7\rm , p. 509}. 
Thus, the lemma follows directly from Faltings' 
\cite{\bf 7\rm , Thm. I.1, p. 509}. 
Though not explicitly stated in Faltings' theorem, 
it follows from its proof that the homomorphism $\Phi_{F,E}$ is also 
additive on $\Cal I$ and $\Cal J$, in a sense that will be left to the 
reader to state.\qed

\vskip0.4cm

Let $S$ be a scheme. Let $\Cal H$ be an $S$-flat coherent sheaf on 
$X\times S$ with constant relative Euler characteristic and multirank 
over $S$. Let $H_0$ be a coherent sheaf on $X$ such that 
$\underline r_{H_0}=\underline r_{\Cal H(s)}$ and 
$\chi(H_0)=\chi(\Cal H(s))$ for every $s\in S$. Let 
$E$ be a vector bundle on $X$. Fix an isomorphism 
$\Cal D(H_0\otimes E)\cong k$. It follows from 
Lemma 1 that, given a vector bundle $F$ on $X$, with $r_F=tr_E$ and 
$\det F\cong(\det E)^{\otimes t}$ for some $t>0$, 
and an isomorphism $\Cal D(H_0\otimes F)\cong k$, 
there is a canonical isomorphism,
$$
\Cal D(\Cal H\otimes F)\cong\Cal D(\Cal H\otimes E)^{\otimes t}.\tag{1.1}
$$
Of course, the isomorphism in (1.1) depends on the choice of $H_0$. 

Assume now that $\chi(E\otimes \Cal H(s))=0$ for 
every $s\in S$. If $F$ is a vector bundle on $X$ with $r_F=tr_E$ and 
$\det F\cong(\det E)^{\otimes t}$ for some $t>0$, 
then we may regard $\theta_F(\Cal H)$ as a section 
of $\Cal L_E(\Cal H)^{\otimes t}$ under the isomorphism in 
(1.1). The induced section 
$\theta_F(\Cal H)\in H^0(S,\Cal L_E(\Cal H)^{\otimes t})$ 
is well defined modulo $k^*$. For every $t>0$, let 
$$
V^t_E(\Cal H)\subseteq H^0(S,\Cal L_E(\Cal H)^{\otimes t})
$$
be the subvectorspace generated by the sections $\theta_F(\Cal H)$. 
We set $V^0_E(\Cal H):=k$. 
It follows from the additive 
properties of the determinant of cohomology and its associated global section, 
and Lemma 1, that, if
$$
0 @>>> F_1 @>>> F @>>> F_2 @>>> 0
$$
is a short exact sequence of vector bundles on $X$, such that 
$r_{F_i}=t_ir_E$ and $\det F_i\cong(\det E)^{\otimes t_i}$ for 
$i=1,2$, then 
$$
\theta_F(\Cal H)=\theta_{F_1}(\Cal H)\otimes\theta_{F_2}(\Cal H)
$$
(modulo $k^*$) inside $H^0(S,\Cal L_E(\Cal H)^{\otimes (t_1+t_2)})$. 
Therefore, the tensor-product multiplication in
$$
\Gamma_E(\Cal H):=\bigoplus_{t\geq 0} H^0(S,\Cal L_E(\Cal H)^{\otimes t})
$$
restricts to a multiplication in
$$
V_E(\Cal H):=\bigoplus_{t\geq 0} V^t_E(\Cal H),
$$
showing that $V_E(\Cal H)$ is a graded $k$-subalgebra of 
$\Gamma_E(\Cal H)$. The ring $V_E(\Cal H)$ is called a 
\sl ring of theta functions\rm .

\vskip0.8cm

\bf 3. The main lemma. \rm Fix a vector bundle $E$ on $X$ of 
rank $r>0$ and multidegree $\underline d$. 
A non-zero torsion-free sheaf $I$ on $X$ is called \sl semistable \rm (resp. 
\sl stable\rm ) \sl with respect to $E$ \rm if: 

(a)\hskip0.2cm $\chi(I\otimes E)=0$;

(b)\hskip0.2cm $\chi(K\otimes E)\leq 0$ (resp. $\chi(K\otimes E)<0$) 
for every proper subsheaf $K\subsetneqq I$.

\parindent=0pt

The notions of stability and 
semistability are numerical, depending on the 
\sl multi-slope \rm $\underline{\mu}:=\underline d/r$ of $E$, rather 
than on $E$ itself. We call $\underline{\mu}$ the \sl polarization \rm. 
We shall fix $E$, and the ensuing polarization, for the remainder of the 
article. 

\parindent=10pt

If $X$ is irreducible, then a 
torsion-free sheaf $I$ on $X$ is semistable (resp. stable) in the usual sense 
if and only if $I$ is semistable (resp. stable) 
with respect to any (hence every) non-zero vector bundle $E$ on $X$ with 
$\chi(I\otimes E)=0$. 

Given any semistable sheaf $I$ on $X$, we can construct a 
\sl Jordan-H\"older filtration\rm ,
$$
0=I_0\subsetneqq I_1 \subsetneqq \dots \subsetneqq I_{q-1} 
\subsetneqq I_q=I,
$$
defined to be a filtration 
where each quotient $J_s:=I_s/I_{s-1}$ is torsion-free and stable, for 
$s=1,\dots, q$. The 
above filtration is not unique, but its graded sheaf,
$$
\text{Gr}(I):=J_1\oplus J_2\oplus\dots\oplus J_q,
$$
is unique by the Jordan-H\"older Theorem.

The following theorem can be thought of as a cohomological characterization of 
semistable sheaves.

\vskip0.4cm

\smc Theorem 2. \rm (Faltings)\hskip0.4cm \sl 
Let $I$ be a non-zero torsion-free sheaf on $X$. Then, 
$I$ is semistable if and only if there is a 
non-zero vector bundle $F$ on $X$ such that:

{\rm (a)}\hskip0.2cm $r\underline d_F=r_F\underline d$;

{\rm (b)}\hskip0.2cm $h^0(X,I\otimes F)=h^1(X,I\otimes F)=0$.

\parindent=0pt

More precisely, if $I$ is semistable, and $\{L_t|\  t>0\}$ is a 
sequence of invertible sheaves on $X$ with 
$\underline d_{L_t}=t\underline d$ for every $t>0$, then 
there is a vector bundle $F$ on $X$ such that, in addition to 
{\rm (a)} and {\rm (b)}:

\parindent=10pt

{\rm (c)}\hskip0.2cm $r_F=tr$ and $\det F\cong L_t$ for a certain $t>0$.

\vskip0.4cm

\sl Proof.\hskip0.4cm \rm As in \cite{\bf 20\rm , Lemma 3.1, p. 166} for 
the ``only if'' part of the first statement, and 
\cite{\bf 20\rm , Lemma 8.3, p. 195} for the ``if'' part. The 
last statement can be proved from the first using 
\cite{\bf 5\rm , Lemma 4} and the argument in the 
proof of \cite{\bf 5\rm , Thm. 5}.\qed

\vskip0.4cm

\sl Remark 3.\hskip0.4cm \rm As 
it could be expected, there is an upper bound $R$ for the rank of $F$ in the 
statements of the above theorem, and this bound depends only on 
the numerical invariants attached to the statements. 
(Similar remark can be made regarding every statement in this article that 
asserts the existence of a vector bundle on $X$ with certain properties.) 
The existence of 
$R$ follows easily from the fact that the family of all 
semistable sheaves on $X$ is bounded. 
In the case where $X$ is non-singular, Le Potier \cite{\bf 13} 
was the first to give an explicit bound $R$. His bound was later 
improved by Hein \cite{\bf 10}. In the general case, one can expect the 
bound to depend on the kind of singularities, as it was the case for 
the compactified Jacobian (cf. \cite{\bf 5}). Anyhow, these general bounds 
seem to be far from sharp (see Beauville's survey article \cite{\bf 3}).

\vskip0.4cm

A $n$-uple of integers 
$\underline{\epsilon}$ is called a \sl deformation of the 
polarization \rm given by $E$ if there is an integer 
$m$ with $m\underline d\equiv \underline{\epsilon} \mod r$. 
Equivalently, $\underline{\epsilon}$ is a deformation of the 
polarization if and only 
if there is a vector bundle $F$ on $X$ such that $r_F\underline d-
r\underline d_F=\underline{\epsilon}$. We say that $\underline{\epsilon}$ 
is \sl non-negative \rm if its components are non-negative. Note that, if 
$\underline{\epsilon}$ is a deformation of the polarization, 
and $I$ is a coherent sheaf 
on $X$ with $\chi(I\otimes E)=0$, 
then $\underline{\epsilon}\cdot\underline r_{I}$ is a multiple of $r$.

Let $S$ be a scheme, and $\Cal F$ an 
$S$-flat coherent sheaf on $X\times S$. We say that $\Cal F$ is 
\sl complete \rm over $S$ if the 
Kodaira-Spencer map,
$$
\delta_s\: T_{S,s} @>>> \text{Ext}^1_X(\Cal F(s),\Cal F(s)),
$$
of $\Cal F$ at $s$ is surjective for every $s\in S$. As 
remarked in the proof of \cite{\bf 7\rm , Thm. I.2, p. 514}, it is easy to 
show that, given any vector bundle $F$ on $X$, there are a connected, 
non-singular scheme $S$ and a vector bundle $\Cal F$ on $X\times S$ such that 
$\Cal F$ is complete over $S$, and 
$\Cal F(s)\cong F$ for a certain $s\in S$.

\vskip0.4cm

\smc Lemma 4.\hskip0.4cm \sl Let 
$I_1,\dots,I_q$ be a finite collection of 
non-isomorphic stable sheaves on $X$. Let 
$\underline{\epsilon}$ be a deformation of the polarization such that 
$\underline{\epsilon}\cdot\underline r_{I_j}\geq 0$ for 
$j=1,\dots,q$. Then, there is a vector bundle $F$ on $X$ such that:

{\rm (a)}\hskip0.2cm $\underline{\epsilon}=r_F\underline d-r\underline d_F$;

{\rm (b)}\hskip0.2cm $H^0(X,I_j\otimes F)=0$ for $j=1,\dots,q$;

{\rm (c)}\hskip0.2cm the natural homomorphism,
$$
\bigoplus_{j=1}^q (I_j\otimes H^1(X,I_j\otimes F)^*) @>>> 
F^*\otimes\omega,
$$
is injective with torsion-free cokernel.

\vskip0.4cm

\sl Proof.\hskip0.4cm \rm Applying 
Thm. 2 to $I_1\oplus\dots\oplus I_q$, we 
obtain a vector bundle $G$ on $X$ such that $r\underline d_G=r_G\underline d$ 
and
$$
h^0(X,I_j\otimes G)=h^1(X,I_j\otimes G)=0\tag{4.1}
$$
for $j=1,\dots,q$.

Let $G_1$ be a vector bundle on $X$ with $r_{G_1}\underline d-r
\underline d_{G_1}=\underline{\epsilon}$. Let $L$ be an ample invertible 
sheaf on $X$ such that $G\otimes L$ and $G_1\otimes L$ are 
generated by global sections. Let $S_1$ be a smooth, connected scheme, 
and $\Cal G_1$ be a vector bundle on $X\times S_1$, complete over $S_1$, 
such that $\Cal G_1(s)\cong G_1$ for a certain $s\in S_1$. Replacing 
$S_1$ by an open dense subscheme, if necessary, we may assume that, for 
every $s\in S_1$, the sheaf $\Cal G_1(s)\otimes L$ is 
generated by global sections, and $h^1_j:=h^0(X,I_j\otimes\Cal G_1(s))$ 
is constant for every $j$. Let $G_2:=\Cal G_1(s_1)\oplus G$, for 
any fixed $s_1\in S_1$. Applying the above procedure to $G_2$, in the 
place of $G_1$, and then 
repeatedly, we obtain a sequence, $(S_m,\Cal G_m)$, of smooth, connected 
schemes $S_m$ and vector bundles $\Cal G_m$ on $X\times S_m$ such 
that, for every $m>0$, the sheaf $\Cal G_m$ is complete over $S_m$; 
$$
h^m_j:=h^1(X,I_j\otimes\Cal G_m(s))
$$
does not depend on $s\in S_m$, for 
every $j$; and $\Cal G_m(s)\otimes L$ is generated by global sections 
for each $s\in S_m$. 
Actually, it follows from our construction and (4.1) that 
$h^{m+1}_j\leq h^m_j$ for every $m>0$ and $j=1,\dots,q$. Thus, 
reenumerating 
the sequence $(S_m,\Cal G_m)$, if necessary, we may assume that 
$h_j:=h^m_j$ does not depend on $m>0$. 
As in the proof of \cite{\bf 20\rm , Lemma 3.1, p. 166}, we have 
that the bilinear composition map,
$$
\text{Hom}_X(\Cal G_m^*(s),I_j)
\times\text{Hom}_X(I_j,\Cal G_m^*(s)\otimes\omega) 
@>>> \text{Hom}_X(\Cal G_m^*(s),\Cal G_m^*(s)\otimes\omega),\tag{4.2}
$$
is zero for every $m>0$, every $s\in S_m$ and $j=1,\dots,q$.

For every $m>0$ and every $j$, 
let $\Cal V^m_j:=R^1p_{m*}(I_j\otimes\Cal G_m)$, where we denote by 
$p_m\: X\times S_m @>>> S_m$ the projection map. Since 
$S_m$ is reduced, it follows from the properties of $\Cal G_m$ that 
$\Cal V^m_j$ is locally free of rank $h_j$. (We may actually assume that 
$\Cal V^m_j$ is free.) Consider the canonical homomorphism,
$$
\lambda_m\: \bigoplus_{j=1}^q (I_j\otimes (\Cal V^m_j)^*) @>>> 
\Cal G_m^*\otimes\omega,
$$
for every $m>0$. 
Replacing $S_m$ by an open dense subscheme, if necessary, we may assume that 
the cokernel $\Cal C_m$ of $\lambda_m$ is flat over $S_m$. Let 
$\Cal J_m$ denote the image of $\lambda_m$. Then, 
$\Cal J_m$ is a relatively torsion-free quotient of 
$\oplus (I_j\otimes(\Cal V^m_j)^*)$ over $S_m$. 
In fact, for every 
$s\in S_m$, the sheaf $\Cal J_m(s)$ is the smallest subsheaf of 
$\Cal G_m^*(s)\otimes\omega$ through which all homomorphisms 
$I_j @>>> \Cal G_m^*(s)\otimes\omega$, for all $j$, factor. 
The sheaf $\Cal C_m$ is not necessarily relatively torsion-free over $S_m$, 
but we may assume that $\Cal C_m^{\omega\omega}$ is (see Sect. 1.1). 
Moreover, we may assume that $\Cal C_m^{\omega\omega}(s)$ is the 
largest torsion-free quotient of $\Cal C_m(s)$ for every $s\in S_m$. Let 
$\Cal K_m\subseteq\Cal G_m^*\otimes\omega$ denote the kernel of the 
composition,
$$
\Cal G_m^*\otimes\omega \twoheadrightarrow \Cal C_m \twoheadrightarrow 
\Cal C_m^{\omega\omega},
$$
for every $m>0$. Then, $\Cal K_m$ is a 
relatively torsion-free sheaf on $X\times S_m$ over $S_m$. 
Moreover, it is clear that $\Cal J_m\subseteq\Cal K_m$, with 
$S_m$-flat quotient 
$\Cal K_m/\Cal J_m$ of relative finite length over $S_m$. 

Let $s\in S_m$. Since $\Cal J_m(s)$ is a quotient of 
$$
I:=\bigoplus_{j=1}^q I_j^{\oplus h_j},
$$ 
then there is a lower bound $a$, independent of $m$ and $s\in S_m$, for 
$\chi(\Cal K_m(s))$. On the other hand, since $\Cal K_m(s)
\subseteq\Cal G_m^*(s)\otimes\omega$ with torsion-free quotient, and 
$\Cal G_m(s)\otimes L$ is generated by global sections, 
then there is an upper bound $A$, independent of $m$ and $s\in S_m$, for 
$\chi(\Cal K_m(s))$. To conclude, 
there are finitely many Hilbert polynomials that 
$\Cal K_m(s)$ can have. More precisely, the 
Hilbert polynomial $P_m(T)$ of $\Cal K_m(s)$ with respect to $L$ is of 
the form:
$$
P(T)=d+T\underline d_L\cdot\underline r, \text{ where } 
a\leq d\leq A \text{ and } \underline 0\leq\underline r\leq\underline 
r_I.
\tag{4.3}
$$

Let $Q\subseteq\text{Quot}_{I^{\omega}}$ denote the subscheme of 
Grothendieck's Quot-scheme which parametrizes quotients of $I^{\omega}$ with 
Hilbert polynomial $P(-T)-P_I(-T)$ with respect to $L$, 
where $P(T)$ ranges through the polynomials of the form (4.3). 
Let $e:=\dim Q$. Note that $e$ depends only on $I$ 
and $L$. Since $\Cal V^m_j$ is free, choosing a 
basis of $\Cal V^m_j$ for every $j$ 
we get that the induced embedding $\Cal K^{\omega}_m
\hookrightarrow I^{\omega}\otimes\Cal O_{S_m}$ 
defines a morphism $g_m\: S_m @>>> Q$. 
As in the proof of \cite{\bf 20\rm , Lemma 3.1, p. 166}, 
since $\Cal G_m$ is complete over $S_m$, we have 
that the image of the bilinear 
composition map,
$$
\text{Hom}_X(\Cal G_m^*(s),\Cal K_m(s))
\times\text{Hom}_X(\Cal K_m(s),\Cal G_m^*(s)\otimes\omega) 
@>>> \text{Hom}_X(\Cal G_m^*(s),\Cal G_m^*(s)\otimes\omega),
$$
is contained in a subspace of dimension at most $e$, for every $s\in S_m$. 
Therefore, since 
$\Cal K_m(s)\subseteq \Cal G_m^*(s)\otimes\omega$, then 
$h^0(X,\Cal K_m(s)\otimes\Cal G_m(s))\leq e$, and hence:
$$
\chi(\Cal K_m(s)\otimes\Cal G_m(s))\leq e\tag{4.4}
$$
for every $m>0$ and $s\in S_m$. 
On the other hand, it follows from our construction of $(S_m,\Cal G_m)$ that
$$
\chi(\Cal K_m(s)\otimes\Cal G_m(s))=(\text{rk }\Cal G_m/r)
\chi(\Cal K_m(s)\otimes E)-(1/r)\underline{\epsilon}\cdot
\underline r_{\Cal K_m(s)}.\tag{4.5}
$$
for every $m>0$ and $s\in S_m$. Combining (4.4) and (4.5), since 
$\text{rk }\Cal G_m\to\infty$ as $m\to\infty$, we get that 
$\chi(\Cal K_m(s)\otimes E)=0$ for every $s\in S_m$, if $m>>0$. 
Since $\Cal K_m(s)\supseteq \Cal J_m(s)$, 
and $\Cal J_m(s)$ is a quotient of the semistable sheaf $I$, 
then $\Cal K_m(s)=\Cal J_m(s)$, and $\Cal J_m(s)$ is 
semistable for every $s\in S_m$, if $m>>0$. 

Fix $s\in S_m$, for $m>>0$. Let $F:=\Cal G_m(s)$ and $J:=\Cal J(s)$. 
Of course, $F$ meets condition (a) in the statement of the lemma. 
Let (after a choice of bases):
$$
\lambda:=\lambda_m(s)\: I @>>> F^*\otimes\omega.
$$
Since $J:=\text{im}(\mu)$ is semistable, and $I$ is a direct sum of 
stable sheaves, then (after a new choice of bases):
$$
J\cong\bigoplus_{j=1}^q I_j^{\oplus h'_j}
$$
for certain integers $h'_j\leq h_j$, and 
the surjective homomorphism, $\lambda\: I \twoheadrightarrow J$, 
has the form:
$$
\lambda=\bigoplus_{j=1}^q (\text{id}_{I_j}\otimes \phi_j),
$$ 
where $\phi_j\: k^{\oplus h_j} @>>> k^{\oplus h'_j}$ is a 
linear surjective homomorphism, for 
$j=1,\dots,q$. We claim that $\phi_j$ is injective for $j=1,\dots,q$. 
Indeed, if there were 
$\underline a\in k^{\oplus h_j}$ 
such that $\phi_j(\underline a)=0$, then the zero homomorphism 
would be written as a linear combination 
with coefficients $a_1,\dots,a_{h_j}$ of the homomorphisms forming 
a basis of $\text{Hom}_X(I_j,F^*\otimes\omega)$. Consequently, 
$\underline a=\underline 0$. The upshot is that $\phi$ is an isomorphism, 
thus showing item (c) in the statement of the lemma. 

Since $\lambda$ is injective, 
there is an embedding $I_j \hookrightarrow F^*\otimes\omega$ for 
$j=1,\dots,q$. It follows from the triviality of 
(4.2) that $H^0(X,I_j\otimes F)=0$ for $j=1,\dots,q$. The 
proof of the lemma is complete.\qed

\vskip0.4cm

It follows from Thm. 2 and the proof of Lemma 4 that, if 
$\underline{\epsilon}\equiv 0 \mod r$, then we may choose the vector 
bundle $F$ in the statement of Lemma 4 with $r|r_F$.

\vskip0.8cm

\bf 4.\hskip0.2cm Quasistable sheaves.\hskip0.4cm \rm 
Fix a vector bundle $E$ on $X$, of rank $r>0$ and 
multidegree $\underline d$, and the ensuing polarization. 
Let $I$ be a semistable sheaf on $X$. Let 
$\underline{\epsilon}$ be a deformation of the polarization. We say that 
$I$ is \sl $\underline{\epsilon}$-quasistable 
\rm (with respect to the polarization) 
if:

(a)\hskip0.2cm $\underline{\epsilon}\cdot\underline r_I>0$;

(b)\hskip0.2cm $\underline{\epsilon}\cdot\underline r_J=0$ for 
every proper semistable quotient $J$ of $I$.

\parindent=0pt

We observe that an $\underline{\epsilon}$-quasistable sheaf is 
simple, that is, its automorphisms are homotheties.

\parindent=10pt

It is clear that a stable sheaf $I$ is $\underline{\epsilon}$-quasistable for 
every deformation $\underline{\epsilon}$ of the polarization such that 
$\underline{\epsilon}\cdot\underline r_I>0$. Conversely, 
if a semistable sheaf $I$ is $(r,\dots,r)$-quasistable, then 
$I$ is stable. 

If $X$ is irreducible, then a semistable sheaf is quasistable if and only 
if it is stable. Thus, no new concept is being introduced in this case. 
On the other hand, if $X$ is reducible, then quasistable sheaves do not 
need to be stable. We will see in \cite{\bf 6} that quasistable 
sheaves are useful in providing a fine (with universal sheaf) 
compactification of the (generalized) Jacobian of a reducible curve.

\vskip0.4cm

\smc Proposition 5.\hskip0.4cm \sl Let $I$ be a semistable sheaf on $X$. Let:
$$
0=I_0\subsetneqq I_1\subsetneqq\dots\subsetneqq I_{q-1}\subsetneqq I_q=I
$$
be a Jordan-H\"older filtration of $I$. Let $\underline{\epsilon}$ 
be a deformation of the polarization such that 
$\underline{\epsilon}\cdot\underline r_I>0$. 
Then, the following statements are equivalent:

{\rm (a)}\hskip0.2cm $I$ is $\underline{\epsilon}$-quasistable.

{\rm (b)}\hskip0.2cm $I_s$ is 
$\underline{\epsilon}$-quasistable for $s=1,\dots,q$.

{\rm (c)}\hskip0.2cm $\underline{\epsilon}\cdot\underline r_{I_1}=\dots=
\underline{\epsilon}\cdot\underline r_{I_q}$, and the short exact sequence:
$$
0 @>>> I_{s-1} @>>> I_s @>>> I_s/I_{s-1} @>>> 0
$$
is not split for $s=2,\dots,q$.

\vskip0.4cm

\sl Proof.\hskip0.4cm \rm We 
assume (a) and prove (b). By descending induction, 
it is enough to prove that $I_{q-1}$ is $\underline{\epsilon}$-quasistable. 
In fact, let 
$J$ be a proper semistable quotient of $I_{q-1}$. We must show that 
$\underline{\epsilon}\cdot\underline r_J=0$. Let 
$K:=\text{ker}(I_{q-1} \twoheadrightarrow J)$. Since $I_{q-1}$ and $J$ are 
semistable, then so is $K$. Since $I$ is semistable, then so is 
$J':=I/K$. Of course, we have a natural embedding $J\subseteq J'$. 
Let $J'':=J'/J$. Since 
$K$ is a proper subsheaf of $I$, and $I$ is 
$\underline{\epsilon}$-quasistable, then: 
$$
\underline{\epsilon}\cdot\underline r_{J''}=
\underline{\epsilon}\cdot\underline r_{J'}=0.
$$
Thus, $\underline{\epsilon}\cdot\underline r_J=0$, completing the 
proof of (b).

We assume (b) and prove (c). Since $I_s$ is $\underline{\epsilon}$-quasistable 
for $s=1,\dots,q$, it is clear that
$$
\underline{\epsilon}\cdot\underline r_{I_1}=\dots=
\underline{\epsilon}\cdot\underline r_{I_q}.
$$
In addition, 
if the exact sequence in (c) were split for a certain $s>1$, then we 
would have that $I_{s-1}$ is a proper semistable quotient of $I_s$ with 
$\underline{\epsilon}\cdot\underline r_{I_{s-1}}>0$, contradicting 
the fact that $I_s$ is $\underline{\epsilon}$-quasistable. 
The proof of (c) is complete.

We assume (c) and prove (a). By induction, we may assume that 
$I_{q-1}$ is $\underline{\epsilon}$-quasistable. 
Let $\lambda\: I \twoheadrightarrow J$ be 
a proper surjective homomorphism, where $J$ is semistable. Assume 
by contradiction that $\underline{\epsilon}\cdot\underline r_J\neq 0$. 
Let $J':=\lambda(I_{q-1})$ and $K:=\text{ker}(\lambda)$. Of course, 
$J/J'$ is a quotient of $I/I_{q-1}$. Since $I/I_{q-1}$ is stable, then 
either $J=J'$ or $J/J'\cong I/I_{q-1}$. Either way, since 
$\underline{\epsilon}\cdot\underline r_I=\underline{\epsilon}\cdot
\underline r_{I_{q-1}}$, we have that 
$\underline{\epsilon}\cdot\underline r_{J'}=\underline{\epsilon}\cdot
\underline r_{J}\neq 0$. 
Since $I_{q-1}$ is $\underline{\epsilon}$-quasistable, 
then $J'\cong I_{q-1}$. But then, the sequence
$$
0 @>>> I_{q-1} @>>> I @>>> I/I_{q-1} @>>> 0
$$
splits, contradicting the hypothesis of (c). The proof is complete.\qed

\vskip0.4cm

The following theorem can be thought of as a cohomological characterization 
of $\underline{\epsilon}$-quasistability, for any deformation 
$\underline{\epsilon}$ of the polarization. Therefore, the theorem provides 
a fortiori a characterization of stability.

\vskip0.4cm

\smc Theorem 6.\hskip0.4cm \sl Let $I$ be a semistable sheaf on $X$. Let 
$\underline{\epsilon}$ be a deformation of the polarization such that 
$\underline{\epsilon}\cdot\underline r_I>0$. If $I$ is 
$\underline{\epsilon}$-quasistable, then 
there is a vector bundle $F$ on $X$ such that:

{\rm (a)}\hskip0.2cm $r\underline d_F=r_F\underline d-\underline{\epsilon}$;

{\rm (b)}\hskip0.2cm $H^0(X,I\otimes F)=0$;

{\rm (c)}\hskip0.2cm the canonical homomorphism,
$$
I\otimes H^1(X,I\otimes F)^* @>>> F^*\otimes\omega,
$$
is injective.

The converse is true if $\underline{\epsilon}$ is non-negative.

\vskip0.4cm

\sl Proof.\hskip0.4cm \rm Assume 
that $I$ is $\underline{\epsilon}$-quasistable. 
Let
$$
0=I_0\subsetneqq I_1\subsetneqq\dots\subsetneqq I_{q-1}\subsetneqq I_q=I
\tag{6.1}
$$
be a Jordan-H\"older filtration of $I$. Let $J_s:=I_s/I_{s-1}$ for 
$s=1,\dots,q$. 
Since $I$ is $\underline{\epsilon}$-quasistable, it 
follows from Prop. 5 that 
$\underline{\epsilon}\cdot\underline r_{J_s}=0$ for $s=2,\dots,q$. 
Consequently, $\underline{\epsilon}\cdot\underline r_{J_1}=\underline{\epsilon}
\cdot\underline r_I>0$. It follows from Lemma 4 that there is a 
vector bundle $F$ on $X$ such that 
$r\underline d_F=r_F\underline d-\underline{\epsilon}$;
$$
h^0(X,J_s\otimes F)=h^1(X,J_s\otimes F)=0
$$
for $s=2,\dots,q$;
$$
h^0(X,J_1\otimes F)=0 \text{\  \  and \  \  } 
h^1(X,J_1\otimes F)=\underline{\epsilon}\cdot\underline r_I/r;
$$
and the natural homomorphism,
$$
\lambda_1\:I_1\otimes H^1(X,I_1\otimes F)^* @>>> F^*\otimes\omega,
$$
is injective. It follows now from the long exact sequences 
in cohomology associated with the filtration (6.1) that 
$H^0(X,I\otimes F)=0$. We shall prove by induction 
on $q$ that the canonical homomorphism (denoted henceforth by 
$\lambda$) of (c) is injective. 
The initial induction step ($q=1$) is given by the injectivity 
of $\lambda_1$. We may now assume that the natural 
homomorphism,
$$
\lambda_{q-1}\: I_{q-1}\otimes H^1(X,I_{q-1}\otimes F)^* @>>> F^*\otimes\omega,
$$
is injective. Note first that
$$
H^1(X,I_1\otimes F)=\dots=H^1(X,I_{q-1}\otimes F)=H^1(X,I\otimes F).\tag{6.2}
$$
Choose a basis for $H^1(X,I\otimes F)$. Then, 
we may view $\lambda$ in the form:
$$
\lambda\: I^{\oplus c} @>>> F^*\otimes\omega,
$$
where $c:=\underline{\epsilon}\cdot\underline r_I/r$. 
It follows from (6.2) that 
$\lambda_{q-1}$ is the composition of 
the embedding $I_{q-1}^{\oplus c} \hookrightarrow I^{\oplus c}$ with 
$\lambda$. Hence, by induction hypothesis, $\lambda$ is injective on 
$I_{q-1}^{\oplus c}$. Let $K:=\lambda(I^{\oplus c})$ and 
$K':=\lambda(I_{q-1}^{\oplus c})$. Then, $\lambda$ induces a 
surjective homomorphism, 
$\rho\: J_q^{\oplus c} \twoheadrightarrow K/K'$. Since $J_q$ is stable, 
choosing a different basis for $H^1(X,I\otimes F)$, if necessary, we 
may assume that:
$$
\rho=(\rho_1,\rho_2)\: J_q^{\oplus c_1}\oplus J_q^{\oplus c_2} 
@>>> K/K',
$$
where $\rho_1$ is an isomorphism, and $\rho_2=0$. Thus, 
$\lambda$ is injective on $I^{\oplus c_1}\oplus I_{q-1}^{\oplus c_2}$. 
We need only prove that $c_2=0$. Assume, by contradiction, that $c_2\neq 0$. 
Then, there is a non-zero 
homomorphism $\mu\: I @>>> F^*\otimes\omega$ such that 
$\mu(I)=\mu(I_{q-1})$. Since $\mu$ is injective on $I_{q-1}$, 
then the exact sequence,
$$
0 @>>> I_{q-1} @>>> I @>>> J_q @>>> 0,
$$
splits. We obtain a contradiction with Prop. 5. The 
proof of the first statement of the theorem is complete.

Assume now that there is a vector bundle $F$ on $X$ 
meeting the three conditions in the statement of the theorem. 
Let $J$ be a proper torsion-free quotient of $I$. It follows 
from (c) that $H^1(X,J\otimes F)=0$. 
Hence, $\chi(J\otimes F)\geq 0$. On the other hand, it follows from 
(a) that
$$
r\chi(J\otimes F)=r_F\chi(J\otimes E)-\underline{\epsilon}\cdot\underline r_J.
$$
Combining the last two statements, we obtain that 
$\underline{\epsilon}\cdot\underline r_J\leq r_F\chi(J\otimes E)$. 
Now, if $J$ is semistable, 
then $\chi(J\otimes E)=0$ and, consequently, 
$\underline{\epsilon}\cdot\underline r_J\leq 0$. Assuming that 
$\underline{\epsilon}$ is non-negative, it follows that 
$\underline{\epsilon}\cdot\underline r_J=0$. The proof of the 
theorem is complete.\qed

\vskip0.4cm

Let $\text{\bf U}$ denote the contravariant functor defined by:
$$
\text{\bf U}(S):=\{\text{relatively torsion-free sheaves } \Cal I 
\text{ on } X\times S \text{ over } S\}/\sim
$$
for every scheme $S$, where two relatively torsion-free sheaves $\Cal I$ and 
$\Cal J$ on $X\times S$ over $S$ are considered equivalent if there exists an 
invertible sheaf $N$ on $S$ such that $\Cal I\cong\Cal J\otimes N$. 
Fix an integer $\chi$, and a $n$-uple $\underline m$ such 
that $r\chi+\underline d\cdot\underline m=0$. 
Let $\underline{\epsilon}$ be a deformation of the polarization 
such that $\underline{\epsilon}\cdot\underline m>0$. In this case, 
$\underline{\epsilon}\cdot\underline m\geq r$, with equality 
only if the maximum common divisor of $\chi,m_1,\dots,m_n$ is 1. Let 
$$
\text{\bf U}^{\underline{\epsilon}}(\underline m,\chi)\subseteq 
\text{\bf U}
$$
denote the subfunctor parametrizing $\underline{\epsilon}$-quasistable 
sheaves with Euler characteristic $\chi$ and multirank $\underline m$.

\vskip0.4cm

\smc Proposition 7.\hskip0.4cm \sl The following statements are true:

{\rm (a)}\hskip0.2cm If $\underline{\epsilon}$ is non-negative, then 
the subfunctor 
$\text{\bf U}^{\underline{\epsilon}}(\underline m,\chi)\subseteq 
\text{\bf U}$ is open.

{\rm (b)}\hskip0.2cm If $\text{\bf V}\subseteq 
\text{\bf U}^{\underline{\epsilon}}(\underline m,\chi)$ is an open 
subfunctor of $\text{\bf U}$, and 
$\underline{\epsilon}\cdot\underline m=r$, then 
$\text{\bf V}$ is 
representable by a scheme.

\vskip0.4cm

\sl Proof.\hskip0.4cm \rm The first statement follows easily from 
the fact that properties (a), (b) and (c) in Thm. 6 are open. To show the 
second statement, we need only show that 
$\text{\bf V}$ is 
locally representable by a scheme. Let 
$F$ be a vector bundle on $X$ with 
$r\underline d_F=r_F\underline d-\underline{\epsilon}$. 
Let $\text{\bf V}_F\subseteq 
\text{\bf V}$ denote the open subfunctor parametrizing 
torsion-free sheaves $I$ on $X$ such that $h^0(X,I\otimes F)=0$, and 
the unique (modulo $k^*$) homomorphism, 
$I @>>> F^*\otimes\omega$, is injective. 
It follows from the hypothesis and Thm. 6 
that $\text{\bf V}$ 
is covered by open subfunctors of the form $\text{\bf V}_F$, for 
$F$ running through all vector bundles with 
$r\underline d_F=r_F\underline d-\underline{\epsilon}$. Therefore, 
we need only show that $\text{\bf V}_F$ is representable by a scheme. Let 
$W\subseteq \text{Quot}_{F^*\otimes\omega}$ be the 
open subscheme parametrizing quotients $q\: F^*\otimes\omega 
\twoheadrightarrow C$ such that $\text{ker}(q)$ has Euler characteristic 
$\chi$ and multirank $\underline m$, and is represented in 
$\text{\bf V}_F$. 
Let $\text{\bf W}$ denote the functor of points of $W$. It should be clear 
now how to construct a map of functors 
$\text{\bf W} @>>> \text{\bf V}_F$, and 
show that it is an isomorphism. The proof is complete.\qed

\vskip0.4cm

Let $\text{\bf U}^s(\underline m,\chi)\subseteq 
\text{\bf U}$ denote the subfunctor parametrizing stable 
sheaves with Euler characteristic $\chi$ and multirank $\underline m$.

\vskip0.4cm

\smc Corollary 8.\hskip0.4cm \sl The subfunctor 
$\text{\bf U}^s(\underline m,\chi)\subseteq \text{\bf U}$ is open. 
If the maximum common divisor of $\chi,m_1,\dots,m_n$ 
is 1, then $\text{\bf U}^s(\underline m,\chi)$ is representable by a scheme.

\vskip0.4cm

\sl Proof.\hskip0.4cm \rm For the openness, we apply Prop. 7, 
item (a), with $\underline{\epsilon}=(r,\dots,r)$. 
As for the representability, we chose a $n$-uple of integers 
$\underline{\delta}$, and an integer $s$, such that 
$s\chi+\underline{\delta}\cdot\underline m=-1$. Let 
$\underline{\epsilon}:=s\underline d-r\underline{\delta}$. It is clear 
that $\underline{\epsilon}$ is a deformation of the polarization. 
Moreover, $\underline{\epsilon}\cdot\underline m=r$. Since
$$
\text{\bf U}^s(\underline m,\chi)\subseteq 
\text{\bf U}^{\underline{\epsilon}}(\underline m,\chi),
$$
it follows from Prop. 7, item (b), that 
$\text{\bf U}^s(\underline m,\chi)$ is representable by a scheme. 
The proof is complete.\qed

\vskip0.8cm

\bf 5.\hskip0.2cm Separation lemmas.\hskip0.4cm \rm Fix a vector 
bundle $E$ on $X$, of rank $r>0$ and multidegree $\underline d$, 
and the ensuing polarization.

\vskip0.4cm

\smc Lemma 9.\hskip0.4cm \sl Let $I,J_1,\dots,J_q$ be non-isomorphic 
stable sheaves on $X$. Then, 
there is a vector bundle $F$ on $X$, with $r_F=tr$ and 
$\det F\cong (\det E)^{\otimes t}$ for some $t>0$, such that 
$$
h^0(X,I\otimes F)=h^1(X,I\otimes F)\neq 0,
$$
and
$$
h^0(X,J_s\otimes F)=h^1(X,J_s\otimes F)=0
$$
for $s=1,\dots,q$.

\vskip0.4cm

\sl Proof.\hskip0.4cm \rm Fix $l$ such that $c:=r_l(I)\neq 0$. Let 
$c_s:=r_l(J_s)$ for $s=1,\dots,q$. Let $\underline{\delta}$ 
be the $n$-uple whose unique non-zero component is the $l$-th component, 
with value 1. According to 
Lemma 4 (and the observation thereafter), 
there is a vector bundle $G$ on $X$ such that:

{\rm (a)}\hskip0.2cm $r_G=t_0r$ and 
$\underline d_G=t_0\underline d-\underline{\delta}$, 
for some $t_0>0$;

{\rm (b)}\hskip0.2cm $H^0(X,(I\oplus J_1\oplus\dots\oplus J_q)\otimes G)=0$;

{\rm (c)}\hskip0.2cm the natural homomorphism (after choices of bases),
$$
\lambda\:I^{\oplus c}\oplus J_1^{\oplus c_1}\oplus
\dots\oplus J_q^{\oplus c_q} @>>> G^*\otimes\omega,
$$
is injective.

Let $C:=G^*\otimes\omega/\lambda(I^{\oplus c})$, and 
consider the injective homomorphism,
$$
\mu=(\mu_1,\dots,\mu_q)\: 
J_1^{\oplus c_1}\oplus\dots\oplus J_q^{\oplus c_q} \hookrightarrow C,
$$
induced by $\lambda$. Pick 
a non-singular point $p\in X$, contained in $X_l$, such that 
$\mu(p)$ is injective. If $c_1+\dots+c_q=0$, 
let $\rho\: C(p) \twoheadrightarrow k$ be any linear surjective 
homomorphism. Otherwise, we proceed as follows: 
Let $e_{s,1},\dots, e_{s,c_s}$ be a basis of $J_s(p)$, 
for $s=1,\dots,q$. For every $s=1,\dots,q$, and every 
pair $(j_1,j_2)$ with $1\leq j_1,j_2\leq c_s$, let
$$
v_{s,(j_1,j_2)}:=\mu_s(p)(0,\dots,0,e_{s,j_1},0,\dots,0),
$$ 
where $e_{s,j_1}$ sits in the $j_2$-th position in the above $c_s$-uple. 
Since 
$\mu(p)$ is injective, then the vectors $v_{s,(j_1,j_2)}\in C(p)$ 
are linearly independent. Let $\rho\: C(p) \twoheadrightarrow k$ be a linear 
homomorphism such that $\rho(v_{s,(j_1,j_2)})=0$ if $j_1\neq j_2$, and 
$\rho(v_{s,(j,j)})=1$ for $s=1,\dots,q$ and $j=1,\dots,c_s$.

Put
$$
F_0:=(\text{ker}(G^* \twoheadrightarrow 
G^*(p) \twoheadrightarrow C(p) @>\rho >> k))^*.
$$
(We chose implicitly a trivialization of $\omega$ at $p$. 
The particular choice of trivialization is irrelevant.) 
Since $p\in X$ is non-singular, then $F_0$ is a
vector bundle on $X$ with the same rank as $G$, and 
$\det F_0\cong\det G\otimes\Cal O_X(p)$. In particular, 
$\underline d_{F_0}=t_0\underline d$. 
By construction, every homomorphism $I @>>> G^*\otimes\omega$ 
factors through $F_0^*\otimes\omega$. Hence,
$$
h^0(X,I\otimes F_0)=h^1(X,I\otimes F_0)=r_l(I)\neq 0.
$$
On the other hand, it follows from our 
choice of $\rho$ that all non-zero homomorphisms 
$J_s @>>> G^*\otimes\omega$, for $s=1,\dots,q$, do not factor 
through $F_0^*\otimes\omega$. Hence, 
$$
h^0(X,J_s\otimes F_0)=h^1(X,J_s\otimes F_0)=0
$$
for $s=1,\dots,q$.

According to Thm. 2, there is a vector bundle $F_1$ on $X$, with 
$r_{F_1}=t_1r$ and 
$\det F_1\cong (\det E)^{\otimes t_1+t_0}\otimes (\det F_0)^{-1}$ for 
some integer $t_1>0$, such that
$$
h^0(X,(I\oplus J_1\oplus\dots\oplus J_q)\otimes F_1)=
h^1(X,(I\oplus J_1\oplus\dots\oplus J_q)\otimes F_1)=0.
$$
It is clear that $F:=F_0\oplus F_1$ meets 
the requirements of the lemma. The proof is complete.\qed

\vskip0.4cm

Let $A:=k[\varepsilon]/\varepsilon^2$. Let 
$S:=\text{\rm Spec }(A)$, and let $s$ denote the unique point of $S$. 
Given a coherent sheaf $H$ on $X$, we say that an $S$-flat 
coherent sheaf $\Cal H$ on $X\times S$ is a \sl deformation \rm of $H$ if 
$\Cal H(s)\cong H$. We say that a deformation $\Cal H$ of $H$ is 
\sl trivial \rm if $\Cal H\cong H\otimes\Cal O_S$.

\vskip0.4cm

\smc Lemma 10.\hskip0.4cm \sl Let $\Cal H$ be a deformation of a coherent 
sheaf $H$ on $X$. Let $F$ be a vector bundle on $X$ such that 
$\chi(H\otimes F)=0$. Then, $\Theta_F(\Cal H)=\{s\}$ scheme-theoretically 
if and only if the two conditions below are verified:

{\rm (a)}\hskip0.2cm $h^0(X, H\otimes F)=h^1(X, H\otimes F)=1$;

{\rm (b)}\hskip0.2cm the unique (modulo $k^*$) non-zero homomorphism 
$H @>>> F^*\otimes\omega$ does not extend to a homomorphism 
$\Cal H @>>> F^*\otimes\omega\otimes\Cal O_S$ over $S$.

\vskip0.4cm

\sl Proof.\hskip0.4cm \rm Left to the reader.\qed

\vskip0.4cm

\smc Lemma 11.\hskip0.4cm \sl Let 
$I$ be a stable sheaf on $X$. If $\Cal I$ is a 
non-trivial deformation of $I$, then there is a vector bundle $F$ on $X$, 
with $r_F=tr$ and $\det F\cong (\det E)^{\otimes t}$ for some $t>0$, 
such that $\Theta_F(\Cal I)=\{s\}$ scheme-theoretically.

\vskip0.4cm

\sl Proof.\hskip0.4cm \rm Let $c_i:=r_i(I)$ for 
$i=1,\dots,n$. Let $c:=c_1+\dots+c_n$. 
According to Lemma 4 (and the observation thereafter), there is a 
vector bundle $G$ on $X$ such that:

{\rm (a)}\hskip0.2cm $r_G=t_0r$ and 
$\underline d_G=t_0\underline d-(1,\dots,1)$, for 
some $t_0>0$;

{\rm (b)}\hskip0.2cm $H^0(X,I\otimes G)=0$;

{\rm (c)}\hskip0.2cm the natural homomorphism (after a choice of basis),
$$
\lambda\: I^{\oplus c} @>>> G^*\otimes\omega,
$$
is injective with torsion-free cokernel.

We will express $\lambda$ in the form:
$$
\lambda=(\lambda^1_1,\dots,\lambda^1_{c_1},\lambda^2_1,\dots,\lambda^2_{c_2},
\dots,\lambda^n_1,\dots,\lambda^n_{c_n}).
$$
Since $H^0(X,I\otimes G)=0$, then $\lambda$ extends to a homomorphism,
$$
\tilde{\lambda}\: \Cal I^{\oplus c} @>>> 
G^*\otimes\omega\otimes\Cal O_S,
$$
with $S$-flat cokernel $\Cal C$. Let $C:=\Cal C(s)$, and consider 
the ``tangential'' homomorphism,
$$
\nu=(\nu^1_1,\dots,\nu^n_{c_n})\: I^{\oplus c} @>>> C,
$$
induced by $\tilde{\lambda}$. Since 
the deformation $\Cal I$ is non-trivial, then $\nu\neq 0$. Since $C$ is 
torsion-free, then there is a non-singular point $p\in X$ where 
$\nu(p)\neq 0$. We may assume that $\nu^1_1(p)\neq 0$. 
We may also assume that $p\in X_1$. Of course, we then have that $c_1>0$. 
Let $p_1:=p$, and pick non-singular points $p_2,\dots,p_n\in X$ 
such that $p_i\in X_i$ for $i=2,\dots,n$. Fix trivializations of 
$\omega$ at the points $p_1,\dots,p_n$. For 
$i=1,\dots,n$, let $e^i_1,\dots,e^i_{c_i}$ be a basis of $I(p_i)$. 
We may assume that $\nu^1_1(p_1)(e^1_1)\neq 0$. 
Let $v\in G^*(p_1)$ such that its image in $C(p_1)$ is 
$\nu^1_1(p_1)(e^1_1)$. 
For $i_1,i_2=1,\dots,n$, $j_1=1,\dots,c_{i_1}$ and $j_2=1,\dots,c_{i_2}$, 
let $v^{i_1,i_2}_{j_1,j_2}:=\lambda^{i_2}_{j_2}(p_{i_1})(e^{i_1}_{j_1})$. 
For each fixed $i$, 
the vectors $v^{i,i_2}_{j_1,j_2}\in G^*(p_i)$ are linearly independent. 
Moreover, since $v$ is not 
in the image of $\lambda(p_1)$, then $v$ is not a linear combination of the 
$v^{1,i_2}_{j_1,j_2}\in G^*(p_1)$. For $i=1,\dots,n$, let 
$\rho_i\: G^*(p_i) \twoheadrightarrow 
k$ be a linear surjective homomorphism such that
$$
\rho_i(v^{i,i_2}_{j_1,j_2})=\left\{\aligned 
1,& \text{\  if } i_2=i \text{ and } j_1=j_2,\\
0,& \text{\  otherwise},
\endaligned\right.
$$
with the unique exception that $\rho_1(v^{1,1}_{1,1})=0$. Since $v$ 
is not a linear combination of the $v^{1,i_2}_{j_1,j_2}\in G^*(p_1)$, then 
we may also assume that $\rho_1(v)=1$.

Put
$$
F_0:=(\text{ker}(G^* \twoheadrightarrow \bigoplus_{i=1}^n G^*(p_i) 
@>(\rho_1,\dots,\rho_n)>> k^{\oplus n}))^*.
$$
Since $p_1,\dots,p_n$ are non-singular, then $F_0$ is a vector bundle on $X$ 
with the same 
rank as $G$, and
$$
\det F_0\cong\det G\otimes\Cal O_X(p_1+\dots+p_n).
$$ 
In particular, $\underline d_{F_0}=t_0\underline d$. By construction, 
$h^1(X,I\otimes F_0)=1$, where the unique (modulo $k^*$) homomorphism, 
$\mu^1_1\: I @>>> F_0^*\otimes\omega$, is 
the factorization of $\lambda^1_1$. Let 
$C':=\text{coker}(\mu^1_1)$. Since 
$\rho_1(v^{1,i_2}_{1,j_2})=0$ for $i_2=1,\dots,n$ and $j_2=1,\dots,c_{i_2}$, 
and $\rho_1(v)=1$, then $\nu^1_1$ is not in the subspace of 
$\text{Hom}_X(I,C)$ spanned by $\text{Hom}_X(I,C')$ and 
$\text{Hom}_X(I,G^*\otimes\omega)$. Therefore, $\mu^1_1$ does not 
extend to a homomorphism $\Cal I @>>> F_0^*\otimes\omega\otimes\Cal O_S$ 
over $S$. Applying Lemma 9, we obtain 
that $\Theta_{F_0}(\Cal I)=\{s\}$ scheme-theoretically. 

According 
to Thm. 2, there is a vector bundle $F_1$ on $X$, with $r_{F_1}=t_1r$ and 
$\det F_1\cong(\det E)^{\otimes t_1+t_0}\otimes(\det F_0)^{-1}$ for 
some $t_1>0$, such that
$$
h^0(X,I\otimes F)=h^1(X,I\otimes F)=0.
$$
It is clear that $F:=F_0\oplus F_1$ meets the requirements of the lemma. 
The proof is complete.\qed
 
\vskip0.8cm

\bf 6.\hskip0.2cm The moduli spaces of 
semistable sheaves.\hskip0.4cm \rm Let 
$\underline a:=(a_1,a_2,\dots,a_n)$ be a $n$-uple of positive 
rational numbers with $a_1+a_2+\dots +a_n=1$. 
In \cite{\bf 19\rm , Part 7} Seshadri defined 
a non-zero torsion-free sheaf $I$ on $X$ to be $\underline a$-semistable 
(resp. $\underline a$-stable) if
$$
\chi(K)\leq \frac{\underline a\cdot\underline r_K}
{\underline a\cdot\underline r_I}\chi(I) 
\text{\  \  \  (resp. \  }
\chi(K)< \frac{\underline a\cdot\underline r_K}
{\underline a\cdot\underline r_I}\chi(I)
\text{)} 
$$
for every proper subsheaf $K\subsetneqq I$. 

Let $S(\underline a,\underline m,\chi)$ 
(resp. $S'(\underline a,\underline m,\chi)$) 
denote the set of isomorphism
classes of $\underline a$-semistable (resp. $\underline a$-stable) 
torsion-free 
sheaves of Euler characteristic $\chi$ and multirank 
$\underline m$ on $X$. As in Sect. 3, there is a 
Jordan-H\"older filtration for any 
$\underline a$-semistable sheaf $I$, and we shall denote the associated 
graded sheaf by $\text{Gr}_{\underline a}(I)$. We say that 
two $\underline a$-semistable sheaves $I_1,I_2$ are \sl 
$\underline a$-equivalent \rm if 
$\text{Gr}_{\underline a}(I_1)=\text{Gr}_{\underline a}(I_2)$. 

\vskip0.4cm

\sl Remark 12.\hskip0.4cm \rm 
We shall observe that Seshadri's notion of stability is equivalent to ours. 
Let $I$ be a non-zero 
torsion-free sheaf on $X$ with Euler characteristic $\chi$ and 
multirank $\underline m$. Assume that $\chi>0$. 
Let $\underline a$ be a $n$-uple of 
positive rational numbers such that $a_1+\dots +a_n=1$. 
Let $\underline A:=(A_1,\dots, A_n)$ be a $n$-uple of 
positive integers such that
$$
A_ia_j=A_ja_i \text{ \ for \  } 1\leq i,j \leq n;\tag{12.1}
$$
and $A:=\underline A\cdot\underline m$ is a multiple of $\chi$. 
Since $a_1+\dots+a_n=1$, then it follows from (12.1) that
$$
a_i=\frac{A_i}{A_1+\dots+A_n}\tag{12.2}
$$
for $i=1,\dots,n$. For each $i=1,\dots,n$, 
let $x^i_1,\dots,x^i_{A_i}\in X_i$ be 
non-singular points of $X$. Let
$$
\Cal O_X(1):=\Cal O_X(\sum \Sb 1\leq i\leq n\\ 1\leq j\leq A_i \endSb 
x^i_j).
$$
Of course, $\Cal O_X(1)$ is an ample sheaf on $X$. Let
$$
E:=(\Cal O_X^{\oplus (r/t-1)}\oplus\Cal O_X(-1))^{\oplus t}\tag{12.3}
$$
for a certain integer $t>0$, to be specified later, and $r:=tA/\chi$. 
Since $\underline d_E=-t\underline A$, then
$$
\chi(I\otimes E)=\chi r+\underline m\cdot\underline d_E=0.
$$
If $K\subseteq I$ is a torsion-free subsheaf, then
$$
\chi(K\otimes E)=\chi(K)r-t\underline A\cdot\underline r_K=
t((\chi(K)/\chi)\underline A\cdot\underline m-
\underline A\cdot\underline r_K).
$$
It follows that $\chi(K\otimes E)\leq 0$ (resp. 
$\chi(K\otimes E)<0$) if and only if
$$
\chi(K)\leq \frac{\underline A\cdot\underline r_K}
{\underline A\cdot\underline m}\chi 
\text{\  \  (resp. \  } 
\chi(K)<\frac{\underline A\cdot\underline r_K}
{\underline A\cdot\underline m}\chi
\text{)}. 
$$
Using (12.2), we see that $I$ is 
$\underline a$-semistable (resp. $\underline a$-stable) if 
and only if $I$ is semistable (resp. stable) with respect 
to $E$. Note that, if $I$ is $\underline a$-semistable, then 
$\text{Gr}(I)=\text{Gr}_{\underline a}(I)$.

\vskip0.4cm

\smc Theorem 13. \rm (Seshadri)\hskip0.4cm \sl 
There is a coarse moduli space
for $S'(\underline a,\underline m,\chi)$, whose underlying scheme is a
quasi-projective variety denoted by 
$U^s(\underline a,\underline m,\chi)$. Moreover, 
$U^s(\underline a,\underline m,\chi)$ 
has a natural projective compactification, to be denoted by 
$U(\underline a,\underline m,\chi)$. 
The set $U(\underline a,\underline m,\chi)$ is
isomorphic to the quotient of $S(\underline a,\underline m,\chi)$ 
by the $\underline a$-equivalence relation.

\vskip0.4cm

\sl Sketch of Proof.\hskip0.4cm \rm We outline briefly the aspects we will 
need in Seshadri's proof \cite{\bf 19\rm , Thm. 15, p. 155} of the theorem. 

First note that a torsion-free sheaf $I$ on $X$ is 
$\underline a$-semistable, or $\underline a$-stable, if and only if 
$I\otimes L$ is, for every invertible sheaf $L$ on $X$. 
Thus, we may assume that $\chi>>0$. Let
$$
Q:=\text{Quot}^{\underline m,\chi}(\Cal O_X^{\oplus \chi})
$$
denote Grothendieck's scheme of quotients of $\Cal O_X^{\oplus \chi}$ with 
Euler characteristic $\chi$ and multirank $\underline m$. 
Let $R\subseteq Q$ be the open subscheme parametrizing
the quotients $q\: \Cal O_X^{\oplus \chi} \twoheadrightarrow I$ such that 
$I$ is torsion-free, 
and the induced homomorphism $k^{\oplus \chi} @>>> H^0(X,I)$ 
is an isomorphism. Let $R^{ss}$ (resp. $R^s$) denote the
subset of $R$ parametrizing the quotients $q\: \Cal O_X^{\oplus \chi} 
\twoheadrightarrow I$ such that $I$ is $\underline a$-semistable (resp.
$\underline a$-stable). Seshadri constructs 
$U(\underline a,\underline m,\chi)$ 
(resp. $U^s(\underline a,\underline m,\chi)$) as the 
good quotient (resp. geometric quotient) of $R^{ss}$ (resp. $R^s$) under 
the obvious action of $\text{SL}(\chi)$, as we describe below. 

We shall assume the definitions and notations of Rmk. 12 in 
what follows. Let
$$
Z:=\prod_{i=1}^n\prod_{j=1}^{A_i} \text{Grass}(\chi,m_i).
$$
Define a morphism, $\tau\: R @>>> Z$, by mapping a quotient, 
$q\: \Cal O_X^{\oplus \chi} \twoheadrightarrow I$, represented in $R$ to 
$$
(q(x^1_1),\dots,q(x^1_{A_1}),\dots,q(x^n_1),\dots,q(x^n_{A_n}))\in Z,
$$
where $q(x_j^i)\: k^{\oplus \chi} \twoheadrightarrow I(x^i_j)$ 
is the homomorphism 
induced by $q$ on $x_j^i$, for all $i,j$. 
It is clear that $\tau$ is an $\text{SL}(\chi)$-morphism, 
where $\text{SL}(\chi)$ acts on both $R$ and $Z$ in the obvious way. For 
$i=1,\dots,n$, let $V_i$ denote the tautological quotient sheaf of 
rank $m_i$ on $\text{Grass}(\chi,m_i)$. Consider the following 
invertible sheaf on $Z$:
$$
M:=(p^1_1)^*\bigwedge^{m_1} V_1\otimes\dots\otimes 
(p^1_{A_1})^*\bigwedge^{m_1} V_1\otimes\dots\otimes
(p^n_1)^*\bigwedge^{m_n} V_n\otimes\dots\otimes
(p^n_{A_n})^*\bigwedge^{m_n} V_n,
$$
where
$$
\text{id}_Z=(p^1_1,\dots,p^1_{A_1},\dots,p^n_1,\dots,p^n_{A_n})\: Z 
@>>> \prod_{i=1}^n\prod_{j=1}^{A_i} \text{Grass}(\chi,m_i).
$$
The group $\text{SL}(\chi)$ acts linearly on $M$. 
Thus, we may define stable and semistable 
points on $Z$ with respect to this action. Let $Z^{ss}$ (resp. $Z^s$) 
be the open subscheme of $Z$ parametrizing semistable (resp. stable) 
points on $Z$ with respect to the linear action 
of $\text{SL}(\chi)$ on $M$. By 
G.I.T. 
(\cite{\bf 8\rm , Thm. 1.10, p. 38} or \cite{\bf 15\rm , Thm. 3.21, p. 84}), 
there is a projective 
good quotient of $Z^{ss}$ for the action of $\text{SL}(\chi)$. 
In addition, some power of $M$ descends to an ample sheaf on the 
quotient.

If $\chi$ and $A$ are large enough, then Seshadri \cite{\bf 19\rm , 
Thm. 19, p. 158} 
shows that it is possible to choose the points $x^i_j$ in such a way that: 

{\rm (a)}\hskip0.2cm $\tau$ is injective;

{\rm (b)}\hskip0.2cm $R^{ss}=\tau^{-1}(Z^{ss})$, 
hence $R^{ss}\subseteq R$ is open;

{\rm (c)}\hskip0.2cm $R^s=\tau^{-1}(Z^s)$, hence $R^s\subseteq R$ is open;

{\rm (d)}\hskip0.2cm the 
induced morphism, $\tau^{ss}\: R^{ss} @>>> Z^{ss}$, is proper.

\parindent=0pt

It follows from (a) and (b) that $\tau^{ss}$ is finite. Applying a 
result of Ramanathan's 
\cite{\bf 19\rm , Prop. 26, p. 31}, we get that there is a good 
quotient,
$$
\phi\: R^{ss} @>>> U(\underline a,\underline m,\chi),
$$
for the action of $\text{SL}(\chi)$. In addition, there is an open subscheme, 
$$
U^s(\underline a,\underline m,\chi)\subseteq 
U(\underline a,\underline m,\chi),
$$
such that $R^s=\phi^{-1}(U^s(\underline a,\underline m,\chi))$ and 
$\left. \phi\right|_{R^s} \: R^s
@>>> U^s(\underline a,\underline m,\chi)$ is a geometric quotient of 
$R^s$ under the
action of $\text{SL}(\chi)$. The varieties 
$U(\underline a,\underline m,\chi)$ and 
$U^s(\underline a,\underline m,\chi)$
are the ones mentioned in the statement of the theorem. The variety 
$U(\underline a,\underline m,d)$ is projective, and 
$\left. \tau^*(M)^{\otimes t}\right|_{R^{ss}}$ 
descends to an ample sheaf on $U(\underline a,\underline m,\chi)$ for 
some $t>0$. Our sketch of Seshadri's proof is complete.\qed

\vskip0.4cm

\parindent=10pt

We shall retain the definitions and notations of 
Rmk. 12 and Thm. 13 in what follows. 
Let $q\: \Cal O_{X\times R}^{\oplus \chi} 
\twoheadrightarrow \Cal I$ be the restriction of the 
universal quotient on $X\times Q$ 
to $X\times R$. Clearly,
$$
\tau^*(M)=\bigotimes \Sb 1\leq i\leq n\\ 1\leq j\leq A_i 
\endSb \bigwedge^{m_i} M^i_j,
$$
where $M^i_j:=\left. \Cal I\right|_{x^i_j\times R}$ 
is regarded as a sheaf on $R$ 
under the canonical isomorphism $x^i_j\times 
R \cong R$, for all $i,j$. 
We choose $E$ as in (12.3), with an integer $t>0$ such that 
$\left.\tau^*(M)^{\otimes t}\right|_{R^{ss}}$ descends to 
$U(\underline a,\underline m,\chi)$. 
Let $g\: X\times R @>>> R$ denote the projection morphism. 
By definition of $R$, 
the adjoint homomorphism, 
$\Cal O_R^{\oplus \chi} @>>> g_*\Cal I$, induced by $q$ is an 
isomorphism, and $H^1(X,\Cal I(s))=0$ for every $s\in R$. 
Thus, $\Cal D(\Cal I)=\Cal O_R$, where $\Cal D$ denotes the determinant of 
cohomology with respect to $g$. Using the exact sequence
$$
0 @>>> \Cal I\otimes\Cal O_X(-1) @>>> \Cal I @>>> 
\bigoplus \Sb 1\leq i\leq n\\ 1 \leq j\leq A_i \endSb 
\left. \Cal I\right|_{x^i_j\times R} @>>> 0,
$$
and the additive property of $\Cal D$, we get that
$$
\aligned
\Cal D(\Cal I\otimes E)=& \Cal D(\Cal I)^{\otimes t(r-1)}\otimes 
\Cal D(\Cal I\otimes\Cal O_X(-1))^{\otimes t}\\
=& \Cal D(\Cal I)^{\otimes tr}\otimes(\bigotimes 
\Sb 1\leq i\leq n\\ 1\leq j\leq A_i \endSb \bigwedge^{m_i} M^i_j)^{\otimes t}\\
=& \tau^*(M)^{\otimes t}.
\endaligned
$$
Hence,
$$
\Cal L_E(\left.\Cal I\right|_{X\times R^{ss}})\cong
\left.\tau^*(M)^{\otimes t}\right|_{R^{ss}}.
$$
As we observed in the ``proof'' of Thm. 13, the sheaf 
$\Cal L_E(\left.\Cal I\right|_{X\times R^{ss}})$ 
descends to an ample sheaf on $U(\underline a,\underline m,\chi)$, henceforth 
denoted by $\Cal L_E(\underline m,\chi)$.

Fix a torsion-free sheaf $I_0$ on $X$ with $\chi(I_0)=\chi$ and 
$\underline r_{I_0}=\underline m$. 
It can be verified from the functorial properties of theta functions, and 
Lemma 1, that the sections 
$$
\theta_F(\left.\Cal I\right|_{X\times R^{ss}})\in 
H^0(R^{ss},\Cal L_E(\left.\Cal I\right|_{X\times R^{ss}}))^{\otimes t}),
$$ 
corresponding to vector bundles $F$ on $X$ with $r_F=tr$ and 
$\det F\cong(\det E)^{\otimes t}$, are $\text{SL}(\chi)$-invariant. 
In other words,
$$
V^t_E(\left.\Cal I\right|_{X\times R^{ss}})\subseteq 
H^0(R^{ss},\Cal L_E(\left.\Cal I\right|_{X\times R^{ss}})^{\otimes t})
^{\text{SL}(\chi)}
$$
for every $t>0$. Thus, the sections 
$\theta_F(\left.\Cal I\right|_{X\times R^{ss}})$ 
of $\Cal L_E(\left.\Cal I\right|_{X\times R^{ss}})^{\otimes t}$ 
descend to sections of $\Cal L_E(\underline m,\chi)^{\otimes t}$, 
which we shall denote by $\theta_F(\underline m,\chi)$. The sections 
$\theta_F(\underline m,\chi)$ are well defined modulo $k^*$. 
We denote by 
$\Theta_F(\underline m,\chi)$ the zero-scheme of 
$\theta_F(\underline m,\chi)$. Let 
$$
V_E(\underline m,\chi):=\bigoplus_{t\geq 0} V^t_E(\underline m,\chi)\subseteq 
\bigoplus_{t\geq 0} H^0(U(\underline a,\underline m,\chi),
\Cal L_E(\underline m,\chi)^{\otimes t})=:\Gamma_E(\underline m,\chi)
$$
denote the graded $k$-subalgebra generated by the sections 
$\theta_F(\underline m,\chi)$. Let 
$$
U_{\theta}(\underline a,\underline m,\chi):=
\text{Proj}(V_E(\underline m,\chi)).
$$
Since $V_E(\underline m,\chi)\subseteq \Gamma_E(\underline m,\chi)$, and 
$U(\underline a,\underline m,\chi)=\text{Proj}(\Gamma_E(\underline m,\chi))$, 
then we have a rational map,
$$
\pi\: U(\underline a,\underline m,\chi) @>>> 
U_{\theta}(\underline a,\underline m,\chi).
$$

\vskip0.4cm

\sl Remark 14.\hskip0.4cm \rm We keep the definitions and notations used so 
far in this section. Let $I_0$ be a semistable sheaf on $X$ of Euler 
characteristic $\chi$ and multirank $\underline m$. Suppose that there 
is a proper semistable subsheaf $I_1\subsetneqq I_0$. Let $I_2:=I_0/I_1$. 
Let $\chi_i:=\chi(I_i)$ and $\underline m_i:=\underline r_{I_i}$ for 
$i=1,2$. Let 
$$
Q_i:=\text{Quot}^{\underline m_i,\chi_i}(\Cal O_X^{\oplus \chi_i})
$$
for $i=1,2$. Let
$$
\nu\:Q_1\times Q_2 @>>> Q
$$
be the morphism sending a pair of quotients,
$$
([q_1\: \Cal O_X^{\oplus \chi_1} \twoheadrightarrow J_1],
[q_2\: \Cal O_X^{\oplus \chi_2} \twoheadrightarrow J_2])
$$
to
$$
[q\: \Cal O_X^{\oplus \chi} = \Cal O_X^{\oplus \chi_1}\oplus 
\Cal O_X^{\oplus \chi_2} @>(q_1,q_2)>> J_1\oplus J_2].
$$
It is clear that $\nu$ defines a monomorphism of the 
corresponding functors of points. It follows from 
\cite{\bf 9\rm , 8.11.5} that $\nu$ is a closed embedding.

For $i=1,2$, let 
$R_i\subseteq Q_i$ denote the open subscheme parametrizing quotients 
$q_i\: \Cal O_X^{\oplus \chi_i} \twoheadrightarrow J_i$ such that 
$J_i$ is torsion-free, and the induced homomorphism,
$$
k^{\oplus \chi_i} @>>> H^0(X,J_i),
$$
is an isomorphism. Of course, $\nu^{-1}(R)=R_1\times R_2$. Let
$$
\mu:=\left. \nu\right|_{R_1\times R_2} \: R_1\times R_2 @>>> R
$$
denote the induced closed embedding. Let 
$\Cal I_i$ denote the restriction of the universal quotient on $X\times Q_i$ 
to $X\times R_i$, for $i=1,2$. For $i=1,2$, let 
$R_i^{ss}$ be the open subscheme of $R_i$ parametrizing 
quotients $q_i\: \Cal O_X^{\oplus \chi_i} \twoheadrightarrow J_i$ 
such that $J_i$ is 
$\underline a$-semistable. Of course, $\mu^{-1}(R^{ss})=R_1^{ss}\times 
R_2^{ss}$. Denote by
$$
\mu^{ss}:= \left. \mu\right|_{R_1^{ss}\times R_2^{ss}} \: 
R_1^{ss}\times R_2^{ss} @>>> R^{ss}
$$
the induced closed embedding. Since
$$
(\text{id}_X,\mu)^*\Cal I\cong\Cal I_1\oplus\Cal I_2,
$$
then it follows from properties of theta functions, and 
Lemma 1 (see its proof), that
$$
(\mu^{ss})^*\Cal L_E(\left.\Cal I\right|_{X\times R^{ss}})\cong
\Cal L_E(\left.\Cal I_1\right|_{X\times R^{ss}_1})\otimes
\Cal L_E(\left.\Cal I_2\right|_{X\times R^{ss}_2}),\tag{14.1}
$$
and
$$
(\mu^{ss})^*\theta_F(\left.\Cal I\right|_{X\times R^{ss}})=
\theta_F(\left.\Cal I_1\right|_{X\times R_1^{ss}})\otimes 
\theta_F(\left.\Cal I_2\right|_{X\times R_2^{ss}})\tag{14.2}
$$
(modulo $k^*$) under the identification (14.1), 
for every vector bundle $F$ on $X$, with $r_F=tr$ and 
$\det F\cong(\det E)^{\otimes t}$ for some $t>0$. Of course, $\mu$ is an 
$\text{SL}(\chi_1)\times \text{SL}(\chi_2)$-morphism. So, $\mu^{ss}$ induces
a morphism,
$$
\alpha\: U(\underline a,\underline m_1,\chi_1)\times 
U(\underline a,\underline m_2,\chi_2) @>>> U(\underline a,\underline m,\chi).
$$
Note that $[I]$ is contained in the image of $\alpha$. Under 
$\alpha$, there are 
relations similar to (14.1) and (14.2), namely:
$$
\alpha^*\Cal L_E(\underline m,\chi)\cong
\Cal L_E(\underline m_1,\chi_1)\otimes
\Cal L_E(\underline m_2,\chi_2),\tag{14.3}
$$
and
$$
\alpha^*\theta_F(\underline m,\chi)=\theta_F(\underline m_1,\chi_1)
\otimes\theta_F(\underline m_2,\chi_2)
$$
(modulo $k^*$) under (14.3), 
for every vector bundle $F$ on $X$, with $r_F=tr$ and 
$\det F\cong(\det E)^{\otimes t}$ for some $t>0$.

\vskip0.4cm

\smc Theorem 15.\hskip0.4cm \sl The natural rational map,
$$
\pi\: U(\underline a,\underline m,\chi) @>>> 
U_{\theta}(\underline a,\underline m,\chi),
$$
is defined everywhere and bijective.

\vskip0.4cm

\sl Proof.\hskip0.4cm \rm It follows from Thm. 2 that $\pi$ is defined 
everywhere. Of course, $\pi$ is dominant. Since 
$U:=U(\underline a,\underline m,\chi)$ is complete, then 
$\pi$ is surjective. We shall now 
show that $\pi$ is injective. Let $I$ and $J$ be semistable sheaves on 
$X$ with $\text{Gr}(I)\not\cong\text{Gr}(J)$. We need to show that 
$\pi([I])\neq\pi([J])$. We may assume that $I\cong\text{Gr}(I)$ and 
$J\cong\text{Gr}(J)$. Assume first that there is a stable summand of 
$I$ that does not occur as a stable summand of $J$. 
Then, it follows from Lemma 9 that there is a vector bundle $F$ on $X$, with 
$r_F=tr$ and $\det F\cong(\det E)^{\otimes t}$ for some $t>0$, 
such that $\theta_F(\underline r,\chi)([I])=0$, but 
$\theta_F(\underline r,\chi)([J])\neq 0$. Hence, 
$\pi([I])\neq\pi([J])$. Assume now that that $I$ and $J$ have a 
common stable summand $K$. Let $\chi_1:=\chi(K)$ and 
$\underline m_1:=\underline r_K$. Put $\chi_2:=\chi-\chi_1$ and 
$\underline m_2:=\underline m-\underline m_1$. Write 
$I=K\oplus I_2$ and $J=K\oplus J_2$. Let
$U_i:=U(\underline a,\underline m_i,\chi_i)$ for $i=1,2$. 
The points $[I],[J]\in U$ lie in the image of the points 
$([K],[I_2]),([K],[J_2])\in U_1\times U_2$, respectively, under the morphism
$$
\alpha\: U_1\times U_2 @>>> U
$$
constructed in Rmk. 14. 
By an induction argument, we may assume that there is an integer $t>0$, 
and vector bundles $F$, $G$ on $X$, with $r_F=r_G=tr$ and 
$\det F\cong\det G\cong(\det E)^{\otimes t}$, such that:

{\rm (a)}\hskip0.2cm $\theta_G(\underline m,\chi)([I])\neq 0$ and 
$\theta_G(\underline m,\chi)([J])\neq 0$;

{\rm (b)}\hskip0.2cm 
$$
\frac{\theta_F(\underline m_2,\chi_2)}
{\theta_G(\underline m_2,\chi_2)}([I_2])\neq 
\frac{\theta_F(\underline m_2,\chi_2)}
{\theta_G(\underline m_2,\chi_2)}([J_2]).
$$

\parindent=0pt

By modifying $F$, if necessary, we may also assume that 
$\theta_F(\underline m_1,\chi)([K])\neq 0$. It follows from 
Rmk. 14 that
$$
\aligned 
\frac{\theta_F(\underline m,\chi)}
{\theta_G(\underline m,\chi)}([I])= & 
\frac{\theta_F(\underline m_1,\chi_1)}
{\theta_G(\underline m_1,\chi_1)}([K])
\frac{\theta_F(\underline m_2,\chi_2)}
{\theta_G(\underline m_2,\chi_2)}([I_2])\\ 
\neq & \frac{\theta_F(\underline m_1,\chi_1)}
{\theta_G(\underline m_1,\chi_1)}([K])
\frac{\theta_F(\underline m_2,\chi_2)}
{\theta_G(\underline m_2,\chi_2)}([J_2])\\ 
= & \frac{\theta_F(\underline m,\chi)}
{\theta_G(\underline m,\chi)}([J]).
\endaligned
$$
Hence, $\pi([I])\neq\pi([J])$. 
The proof is complete.\qed

\vskip0.4cm

\parindent=10pt

\smc Theorem 16.\hskip0.4cm \sl Assume that the maximum common divisor of 
$\chi,m_1,\dots,m_n$ is 1. Then, the restriction,
$$
\pi^s\:=\left.\pi\right|_{U^s(\underline a,\underline m,\chi)}\: 
U^s(\underline a,\underline m,\chi) @>>> 
U_{\theta}(\underline a,\underline m,\chi),
$$
is an open embedding.

\vskip0.4cm

\sl Proof.\hskip0.4cm \rm Let $U:=U(\underline a,\underline m,\chi)$ and 
$U_{\theta}:=U_{\theta}(\underline a,\underline m,\chi)$. 
Since $U$ is complete, 
then $\pi$ is proper. Since $\pi$ is bijective, 
then $\pi$ is a homeomorphism, and a finite morphism. To show that 
$\pi^s$ is an open embedding, we need only show that the 
homomorphism of tangent spaces,
$$
d\pi_{[I]}\: T_{[I],U} @>>> T_{\pi([I]),U_{\theta}},
$$
is injective for every stable sheaf $I$ represented in $U$. 
Let $v\in T_{[I],U}$ be a non-zero tangent vector. 
We need to 
show that $d\pi_{[I]}(v)\neq 0$. Since 
$U^s(\underline a,\underline m,\chi)$ is representable, by Cor. 8, then 
$v$ is represented by a deformation $\Cal I$ of $I$. Since $v\neq 0$, 
then $\Cal I$ is not the 
trivial deformation of $I$. It follows from Lemma 11 that there is a 
vector bundle $F$ on $X$, with $r_F=tr$ and $\det F\cong(\det E)^{\otimes t}$ 
for some $t>0$, such that $\Theta_F(\Cal I)=\{s\}$ scheme-theoretically. 
Since $\theta_F(\underline m,\chi)\in V^t_E(\underline m,\chi)$, it 
follows that $d\pi_{[I]}(v)\neq 0$. The proof is complete.\qed

\vskip0.4cm

\sl Remark 17.\hskip0.4cm \rm If 
$X$ is irreducible, and the rank $m$ is coprime with 
the Euler characteristic $\chi$, 
then $U(m,\chi)=U^s(m,\chi)$. So, it follows from Thm. 16 
that $\pi$ is an isomorphism in this case.

\vskip0.4cm

The hypothesis in the statement of Thm. 16 was used in order to identify 
a tangent vector of $U_{\theta}(\underline a,\underline m,\chi)$ at a 
point $[I]$, representing a stable sheaf $I$, with a deformation of $I$. 
It might be that we do not need that hypothesis for such identification, 
as it is the case of the next theorem. 

\vskip0.4cm

\smc Theorem 18.\hskip0.4cm \sl Assume that $X$ is irreducible, and 
the characteristic of the ground field $k$ is 0. Then, the restriction,
$$
\pi^s\:=\left.\pi\right|_{U^s(m,\chi)}\: 
U^s(m,\chi) @>>> U_{\theta}(m,\chi),
$$
is an open embedding.

\vskip0.4cm

\sl Proof.\hskip0.4cm \rm Apply the same proof of Thm. 16, now using 
\cite{\bf 12\rm , Thm. 8.14, p. 141} to guarantee that every tangent 
vector of $U^s(m,\chi)$ is represented by a deformation. The 
proof is complete.\qed

\vskip0.4cm

\sl Remark 19.\hskip0.4cm \rm As 
we had mentioned in Rmk. 3, all the results we 
obtained in this article asserting the existence of a vector bundle 
$F$ on $X$ with certain properties could have been obtained within 
a certain range for the rank of $F$, 
depending only on numerical invariants. Thus, replacing 
$V_E(\underline m,\chi)$ by the subalgebra,
$$
V_E^{[t]}(\underline m, \chi)\subseteq V_E(\underline m,\chi),
$$
generated by all theta functions $\theta_F(\underline m,\chi)$ 
associated with vector bundles $F$ on $X$ with rank $r_F\leq tr$, and letting
$$
U_{\theta}^{[t]}(\underline a,\underline m,\chi):=
\text{Proj}(V_E^{[t]}(\underline m,\chi)),
$$
we have that Theorems 15, 16 and 18 hold for the rational map,
$$
\pi^{[t]}\: U(\underline a,\underline m,\chi) @>>> 
U_{\theta}^{[t]}(\underline a,\underline m,\chi),
$$
as long as $t>>0$. In particular, we have that $\pi^{[t]}$ is 
finite and scheme-theoretically surjective for $t>>0$. Since 
$\pi$ has the same properties, and 
$U_{\theta}(\underline a,\underline m,\chi)$ is the 
inverse limit of $U_{\theta}^{[t]}(\underline a,\underline m,\chi)$ 
as $t\to\infty$, then it follows that
$$
U_{\theta}(\underline a,\underline m,\chi)=
U_{\theta}^{[t]}(\underline a,\underline m,\chi)
$$
for every $t>>0$.

\vskip0.4cm

\sl Acknowledgements.\hskip0.4cm \rm I take the 
opportunity to thank Prof. Kleiman for several comments on 
different versions of the present article. I am also grateful to Prof. 
Ramero for telling me about Faltings' construction of the moduli 
space of semistable vector bundles on a smooth curve via theta functions. 
I would like to thank Waseda University, specially Prof. Kaji, Prof. 
Morimoto, Prof. Ohno and Prof. Hara, for the warm hospitality 
extended during the period this work was begun. I am also grateful 
to Prof. Homma for all the help received from him during the said period.

\vskip0.4cm

\eightsmc Instituto de Matem\'atica Pura e Aplicada, 
Estrada D. Castorina 110, 22460-320 Rio de Janeiro RJ, Brazil


\vskip0.8cm

\Refs

\ref \no 1 \by {\smc A. Altman and S. Kleiman} \paper Compactifying the 
Picard scheme \jour Adv. Math. \vol 35 \yr 1980 \pages 50--112
\endref

\ref \no 2 \by {\smc E. Arbarello, M. Cornalba, P. Griffiths, and J. Harris} 
\paper Geometry of algebraic curves, Vol. I \jour Grundlehren der 
mathematischen Wissenschaften \bf 267\rm , Springer-Verlag, New York, 1985
\endref

\ref \no 3 \by {\smc A. Beauville} 
\paper {\rm ``Vector bundles on curves and 
generalized theta functions: recent results and open problems'' in 
\sl Current topics in complex algebraic geometry} 
\jour {Math. Sci. Res. Inst. Publ. \bf 28\rm , Cambridge Univ. Press, 
Cambridge, 1995, 17--33}
\endref

\ref \no 4 \by {\smc L. Caporaso} \paper A compactification of the universal 
Picard variety over the moduli space of stable curves \jour J. Amer. 
Math. Soc. \vol 7 \yr 1994 \pages 589--660
\endref

\ref \no 5 \by {\smc E. Esteves} \paper Very ampleness for Theta on the 
compactified Jacobian \jour available from the e-print service at 
alg-geom\@ eprints.math.duke.edu/9709005, September 1997; 
to appear in Math. Z.
\endref

\ref \no 6 \by {\smc E. Esteves} \paper Compactifying the relative Jacobian 
over families of reduced curves \jour available from the e-print service at 
alg-geom\@ eprints.math.duke.edu, September 1997
\endref

\ref \no 7 \by {\smc G. Faltings} \paper Stable $G$-bundles and
projective connections \jour J. Algebraic Geometry \vol 2
\yr 1993 \pages 507--568
\endref

\ref \no 8 \by {\smc J. Fogarty. and D. Mumford} 
\paper Geometric Invariant Theory \jour {Ergebnisse der Mathematik und 
ihrer Grenzgebiete \bf 34\rm , Springer-Verlag, Berlin Heidelberg, 1982}
\endref

\ref \no 9 \by {\smc A. Grothendieck with J. Dieudonn\'e} 
\paper \'Elem\'ents de G\'eom\'etrie Alg\'ebrique IV-3 
\jour Inst. Hautes \'Etudes Sci. Publ. \vol 28 \yr 1966
\endref

\ref \no 10 \by {\smc G. Hein} \paper On the generalized theta divisor 
\jour Beitr\"age Algebra Geom. \vol 38 \yr 1997 \pages 95--98
\endref

\ref \no 11 \by {\smc F. Knudsen and D. Mumford} 
\paper The projectivity of the moduli 
space of stable curves I: Preliminaries on ``det'' and ``Div'' 
\jour Math. Scand. \vol 39 \yr 1976 \pages 19--55
\endref

\ref \no 12 \by {\smc J. Le Potier} \paper Fibr\'es vectoriels sur les 
courbes alg\'ebriques \jour {Publications Math\'ematiques de l'Universit\'e 
Paris 7 -- Denis Diderot \bf 35\rm , Universit\'e Paris 7 - Denis Diderot, 
U.F.R. de Math\'ematiques, Paris, 1995}
\endref

\ref \no 13 \by {\smc J. Le Potier} 
\paper \rm ``Module des fibr\'es semi-stables et 
fonctions th\^eta'' in \sl Moduli of vector bundles\rm , 
Lecture Notes in Pure and Appl. Math. \bf 179\rm , Dekker, New York, 1996 
\pages 83--101
\endref 

\ref \no 14 \by {\smc D. Mumford} 
\paper {\rm ``Projective invariants of projective 
structures and applications'' in \sl 
Proc. Intern. Cong. Math. Stockholm} \pages 526--530
\endref

\ref \no 15 \by {\smc P. Newstead} \paper Introduction to moduli problems 
and orbit spaces \jour Tata Inst. Lecture Notes, Springer-Verlag, 1978
\endref

\ref \no 16 \by {\smc R. Pandharipande} 
\paper A compactification over $\overline{M}_g$ 
of the universal moduli space of slope-semistable vector bundles 
\jour J. Amer. Math. Soc. \vol 9 \yr 1996 \pages 425--471
\endref

\ref \no 17 \by {\smc C. S. Seshadri} 
\paper Space of unitary vector bundles on a 
compact Riemann surface \jour Ann. of Math. \vol 85 \yr 1967 
\pages 303--336
\endref

\ref \no 18 \by {\smc C. S. Seshadri} 
\paper {\rm ``Theory of moduli'' in \sl Algebraic Geometry} 
\jour {ed. by R. Hartshorne, Proc. Sympos. Pure Math. \bf 29\rm , 
Amer. Math. Soc., Providence, 1975} \pages 263--304
\endref

\ref \no 19 \by {\smc C. S. Seshadri} \paper Fibr\'es vectoriels
sur les courbes alg\'ebriques \jour Ast\'erisque \vol 96 
\yr 1982
\endref

\ref \no 20 \by {\smc C. S. Seshadri} \paper {\rm ``Vector bundles on
curves'' in \sl Linear algebraic groups and their representations} 
\jour {Contemp. Math. \bf 153\rm , Amer. Math. Soc., Providence, 1993}
\pages 163--200
\endref

\ref \no 21 \by {\smc C. Simpson} \paper Moduli of representations of the 
fundamental group of a smooth projective variety I 
\jour Inst. Hautes \'Etudes Sci. Publ. \vol 79 
\yr 1994 \pages 47--129
\endref

\endRefs
\enddocument